\documentclass[12pt,]{article}
\usepackage{amsmath}
\usepackage{amscd}
\usepackage{amssymb}
\usepackage{latexsym}
\newtheorem{theorem}{Theorem}[section]
\newtheorem{lemma}[theorem]{Lemma}
\newtheorem{corollary}[theorem]{Corollary}
\newtheorem{proposition}[theorem]{Proposition}

\newtheorem{remark}[theorem]{Remark}
\newtheorem{definition}[theorem]{Definition}

\newtheorem{slemma}{Lemma}[subsection]        
\newtheorem{sproposition}[slemma]{Proposition} 
\newtheorem{sdefinition}[slemma]{Definition}  
\newenvironment{proof}{{\em Proof.}\ \ \ }{\unskip\nobreak\hfil\penalty50
\hskip1em\hbox{}\nobreak\hfil $\Box$
\parfillskip=0pt \finalhyphendemerits=0 \par\medskip\noindent}
\newcommand{\bgl}{\begin{equation}}         
\newcommand{\egl}{\end{equation}}
\newcommand{\bgln}{\begin{eqnarray}}        
\newcommand{\egln}{\end{eqnarray}}
\newcommand{\bglnoz}{\begin{eqnarray*}}     
\newcommand{\eglnoz}{\end{eqnarray*}}
\newcommand{\btheo}{\begin{theorem}}
\newcommand{\etheo}{\end{theorem}}
\newcommand{\blemma}{\begin{lemma}}
\newcommand{\elemma}{\end{lemma}}
\newcommand{\bproof}{\begin{proof}}
\newcommand{\eproof}{\end{proof}}
\newcommand{\bbew}{\begin{beweis}}
\newcommand{\ebew}{\end{beweis}}
\newcommand{\bremark}{\begin{remark}\em}
\newcommand{\eremark}{\end{remark}}
\newcommand{\bdefin}{\begin{definition}}
\newcommand{\edefin}{\end{definition}}
\newcommand{\bprop}{\begin{proposition}}
\newcommand{\eprop}{\end{proposition}}
\newcommand{\bcor}{\begin{corollary}}
\newcommand{\ecor}{\end{corollary}}
\newcommand{\mn}{\par\medskip\noindent}

\newcommand{\mc}{\mathcal}

\newcommand{\cC}{\mbox{$\cal C$}}

\newcommand{\lori}{\longrightarrow}
\newcommand{\lole}{\longleftarrow}

\def\SEMI{\mbox{$\times\kern-2pt\vrule height5pt width.6pt \kern3pt $}}

\newcommand{\End}{{\rm End}\,}
\newcommand{\Ker}{{\rm Ker\,}}

\newcommand{\Sp}{{\rm Sp\,}}

\newcommand{\id}{{\rm id}}

\def\Cz{\mathbb{C}}

\def\Nz{\mathbb{N}}

\def\Rz{\mathbb{R}}
\def\Zz{\mathbb{Z}}

\newcommand{\abs}[1]{\lvert#1\rvert}     
\newcommand{\defeq}{\mathrel{:=}}     
\newcommand{\hot}{\mathbin{\hat{\otimes}}}

\begin{document}
\title{Bivariant $K$-theory and the Weyl algebra}%
\author{Joachim Cuntz \\ Mathematisches Institut \\ Universit\"at
M\"unster
\\ Einsteinstr. 62 \\ 48149 M\"unster \\ cuntz@math.
uni-muenster.de}%
\thanks{Research supported by the Deutsche Forschungsgemeinschaft}
\maketitle
\begin{abstract}\noindent
We introduce a new version $kk^{\rm alg}$ of bivariant
$K$-theory that is defined on the category of all locally
convex algebras. A motivating example is the Weyl algebra
$W$, i.e. the algebra generated by two elements satisfying
the Heisenberg commutation relation, with the fine locally
convex topology. We determine its $kk^{\rm alg}$-invariants
using a natural extension for $W$. Using analogous methods the
$kk^{\rm alg}$-invariants can be determined for many other
algebras of similar type.
\end{abstract}
\section{Introduction}
In this article we develop a bivariant topological $K$-theory
for arbitrary locally convex algebras. This covers a very
wide class of examples and even the case of algebras over the
complex numbers without a priori topology. In fact, every
complex algebra $A$ with a countable basis can be considered
as a locally convex algebra with the fine topology defined by
the collection of \underline{all} seminorms on $A$. \\The
construction of the bivariant theory that we present here is
fundamentally the same as the one in \cite{CuDoc} where we
had constructed such a theory for the category of locally
convex algebras with submultiplicative seminorms. As before
it is based on extensions of higher length and ``classifying
maps" for such extensions. We use the opportunity however to
streamline the approach in \cite{CuDoc}. Motivated by the
thesis of A.Thom \cite{Thom} we define the theory as
noncommutative stable homotopy, i.e. as an inductive limit
over suspensions of both variables (using a noncommutative
suspension for the first variable), rather than an inductive
limit over inverse Bott maps as in \cite{CuDoc}. This
simplifies some of the arguments and also clarifies the fact
that Bott periodicity becomes automatic once we stabilize by
(smooth) compact operators in the second variable. As pointed
out in \cite{Thom}, a natural framework
for these constructions is the one of triangulated categories.\\
Just as for the theory $kk$ in \cite{CuDoc} we obtain a
bivariant Chern-Connes character from the new bivariant
$K$-theory $kk^{\rm alg}(A,B)$ to bivariant periodic cyclic
homology $HP_*(A,B)$. In contrast to \cite{CuDoc} however we
loose part of the control over the coefficients $kk^{\rm
alg}(\Cz,\Cz)$ of the theory. Recall that for the theory $kk$
in \cite{CuDoc} we had $kk(\Cz,\Cz)=\Zz$. \\For the new
theory the coefficient ring $R=kk^{\rm alg}(\Cz,\Cz)$ is a
unital ring admitting a unital homomorphism into the field of
complex numbers (via the identification $\Cz=HP_*(\Cz,\Cz)$).
In particular, $R\neq 0$. We leave the determination of $R$
open for the time being.\footnote{Added at the time of
publication: We succeeded recently in showing that $R=\Zz$, at
least if we stabilize by the $C^*$-algebra of compact
operators or by the algebra of Hilbert-Schmidt operators in
place of the algebra $\mathcal K$ of smooth compact
operators. Details will appear elsewhere.} All the abelian groups $kk^{\rm alg}(A,B)$ are
modules over $R$ and the product in $kk^{\rm alg}$ is
$R$-linear. $kk^{\rm alg}$ has the usual properties which
allow to compute $kk^{\rm alg}(A,B)$ in many cases
as $R$-modules.\\
A main motivating example for our theory is the Weyl algebra
$W$ which is the unital complex algebra generated by two
elements $x,y$ satisfying the commutation relation $xy-yx=1$.
It is \underline{the} example of an algebra which cannot be
topologized using submultiplicative seminorms. It is
easy to see that the only submultiplicative seminorm on $W$
is 0. We construct an extension $W'$ of $W$ which can be used
to compute the $kk^{\rm alg}$-invariants for $W$ (as
$R$-modules). The result is that $W$ is isomorphic to $\Cz$
in the category $kk^{\rm alg}$. Along the same lines $kk^{\rm
alg}$-isomorphisms can be determined for many other algebras
of similar type.\\
Generally speaking, the information on a given locally convex
algebra obtained through $kk^{\rm alg}$ is similar to the one
obtained from periodic cyclic homology. In the first case, we
obtain $R$-modules, in the second case vector spaces over
$\Cz$. Thus torsion or divisibility phenomena can not be
controlled in both cases. Therefore, in the last section, we
describe another bivariant $K$-theory $kk$ which is defined
on the category of locally convex algebras but which
restricts to the theory $kk$ of \cite{CuDoc} on the category
of $m$-algebras. In particular, we get for the coefficients
of this theory $kk(\Cz,\Cz)=\Zz$ and the invariants given by
the theory are $\Zz$-modules and therefore contain more
information. The drawback of the theory is the fact that it
admits long exact sequences only for certain extensions and
that, in general, all computations in this setting become
rather technical and complicated. However, the fundamental
extension that we constructed for the Weyl algebra fits into
the framework of this theory, so that we can compute the
$kk$-invariants for the Weyl algebra at least partially. We
give some examples of applications to the structure of the
Weyl algebra that could not be obtained from periodic cyclic homology.
\section{Locally convex algebras}\label{lconv}
Recall the definition of the projective tensor product by
Grothendieck, \cite{GrTens}, \cite{T}. For two locally convex
vector spaces $V$ and $W$ the projective topology on the
tensor product $V\!\otimes\!W$ is given by the family of
seminorms $p\otimes q$, where $p$ is a continuous seminorm on
$V$ and $q$ a continuous seminorm on $W$. Here $p\otimes q$
is defined by
$$ p\!\otimes\!q\,
(z)=\inf\,\Big\{\mathop{\sum}\limits_{i=1}^{n}p(a_i)q(b_i)\;
\Big|\;\,z=
\mathop{\sum}\limits_{i=1}^n a_i\otimes b_i, a_i\in V, b_i\in
W\Big\}
$$ for $z\in V\otimes W$.
We denote by $V\hat{\otimes}W$ the completion of
$V\!\otimes\!W$ with respect to this family of seminorms. \\
By a locally convex algebra we mean an algebra over $\Cz$
equipped with a complete locally convex topology such
that the multiplication $A\times A\to A$ is (jointly) continuous.\\
This means that, for every continuous seminorm $\alpha$ on
$A$, there is another continuous seminorm $\alpha'$ such that
$$\alpha(xy) \leq \alpha'(x)\alpha'(y)$$
for all $x,y \in A$. Equivalently, the multiplication map
induces a continuous linear map $A \hot A \to A$ from the
completed projective tensor product $A\hot A$.\\ It is
convenient to group the continuous seminorms on a locally
convex algebra into sequences of the form $(\alpha_1,
\alpha_2, \ldots )$ such that $\alpha_i(xy) \leq
\alpha_{i+1}(x)\alpha_{i+1}(y)$ for all $i$. We call such a
sequence a submultiplicative sequence of seminorms. Thus, a
locally convex algebra is a complex algebra with a locally
convex structure determined by a family of submultiplicative
sequences of seminorms.\\ The class of locally convex
algebras is very vast. For instance it contains all algebras
over $\Cz$ generated algebraically by a countable family of
generators if we equip them with the ``fine'' locally convex
topology. The fine topology on a complex vector space $V$ is
given by the family of \underline{all} seminorms on $V$. For
completeness we prove the following proposition, cf. also
\cite{Bour}, chap. II, \S 2, Exercise 5 (I am indebted to
J.Dixmier for this reference). \bprop Every algebra $A$ over
$\Cz$ with a countable basis is a locally convex algebra with
the fine topology.\eprop \bproof It is easy to see that $A$ is
complete with the fine topology. Let $(e_i)_{i\in \Nz}$ be a
basis of $A$ and $\alpha$ a seminorm on $A$. Put
$b_{ij}=\alpha (e_i e_j)$. Then, given $x=\sum \lambda_ie_i$
and $y=\sum \mu_je_j$, we have
$$ \alpha (xy) \leq \sum \abs{\lambda_i\mu_j}b_{ij}$$
Define inductively $\beta_1 = \max \{\sqrt{b_{11}},\,1\}$ and
$$\beta_j = \max\;\Big(
\left\{\frac{b_{kj}}{\beta_k},\,\frac{b_{jk}}{\beta_k}\,
\mathop{\Big|}\,k<j\right\}\cup\{\,\sqrt{b_{jj}},\,
1\}\Big)$$ Then $b_{kj}\leq\beta_j\beta_k$ and
$b_{jk}\leq\beta_j\beta_k$ for all $k\leq j$.\\
Therefore, defining a seminorm $\beta$ on $A$ by
$$\beta (\sum \lambda_ie_i )= \sum \abs{\lambda_i}\beta_i$$
we get that $\alpha (xy) \leq \sum \abs{\lambda_i\mu_j}
\beta_i\beta_j = \beta(x)\beta(y)$ for all $x,y$ in
$A$.\eproof For locally convex algebras $A$ and $B\, $, the
projective tensor product $A \hat{\otimes}B $ is again a
locally
convex algebra.\\
As in \cite{CuDoc} we will use the name $m$-algebra for a
locally convex algebra $A$ where the topology is given by a
family of submultiplicative seminorms, i.e. seminorms
$\alpha$ such that
$$\alpha(xy) \leq \alpha(x)\alpha(y)$$
for all $x,y \in A$. The projective tensor product of two
$m$-algebras is again an $m$-algebra. Most of the standard
algebras that we use in our constructions are $m$-algebras. In
particular we will use the algebras described in the
following subsections.
\subsection{Algebras of differentiable functions and
diffotopy.}\label{df} Let $[a,b]$ be an interval in $\Rz$. We
denote by $\Cz[a,b]$ the algebra of complex-valued
$\cC^\infty$-functions $f$ on $[a,b]$, all of whose
derivatives vanish in $a$ and in $b$ (while $f$ itself may
take arbitrary values in $a$ and $b$). \mn Also the
subalgebras $\Cz(a,b], \Cz[a,b)$ and $\Cz(a,b)$ of
$\Cz[a,b]$, which, by definition consist of functions $f$,
that vanish in $a$, in $b$, or in $a$ and $b$, respectively,
will play an important role. \mn The topology on these
algebras is the usual Fr\'echet topology, which is defined by
the following family of submultiplicative norms $p_n$:
$$
p_n(f)=\|f\|+\|f'\|+ {\textstyle\frac{1}{2}} \| f''\|+\dots+
{\textstyle\frac{1}{n!}} \|f^{(n)}\|
$$
Here of course $\|g\|=\sup \{ |g(t)|\,  \big|\,t\in [a,b]\}$
and $f^{(n)}$ denotes the $n$-th derivative of $f$. \mn We
note that $\Cz[a,b]$ is nuclear in the sense of Grothendieck
\cite{GrTens} and that, for any complete locally convex space
$V$, the space $\Cz[a,b]\hat{\otimes}V$ is isomorphic to the
space of $\cC^\infty$-functions on $[a,b]$ with values in
$V$, whose derivatives vanish in both endpoints, \cite{T},
\S\, 51. Exactly the same comments apply to the algebra
$\cC^\infty (M)$ of smooth functions on a compact smooth
manifold where an $m$-algebra topology is given by seminorms
which are defined locally like the $p_n$ above. \mn Given a
locally convex algebra $A $, we write $A[a,b]$, $A [a,b)$ and
$ A (a,b)$ for the locally convex algebras $A \hat{\otimes}
\Cz[a,b]$, $A \hat{\otimes}\Cz[a,b)$ and
$A\hat{\otimes}\Cz(a,b)$ (their elements are $A$ - valued
$\mathcal C ^\infty$-functions whose derivatives vanish at
the endpoints).
\begin{sdefinition} Two continuous homomorphisms $\alpha,
\beta: A\to B$ between locally convex algebras are called
differentiably homotopic, or {\em diffotopic}, if there is a
family $\varphi_t:A\to B$, $t\in [0,1]$, of continuous
homomorphisms, such that $\varphi_0=\alpha, \varphi_1=\beta$
and such that the map $t\mapsto \varphi_t(x)$ is infinitely
often differentiable for each $x\in A$. An equivalent
condition is that there should be a continuous homomorphism
$\varphi:A\to \cC^\infty ([0,1])\hat{\otimes} B$ such that
$\varphi(x)(0)=\alpha(x), \varphi(x)(1)=\beta(x)$ for each
$x\in A$. \end{sdefinition} Let $h:[0,1]\to [0,1]$ be a
monotone and bijective $\cC^\infty$-map, whose restriction to
$(0,1)$ gives a diffeomorphism $(0,1)\to (0,1)$ and whose
derivatives in $0$ and $1$ all vanish. Replacing $\varphi_t$
by $\psi_t=\varphi_{h(t)}$ one sees that $\alpha$ and $\beta$
are diffotopic if and only if there is a continuous
homomorphism $\psi:A\to \Cz[0,1]\hat {\otimes} B$ such that
$\psi(x)(0)=\alpha(x), \psi(x)(1)=\beta(x)$, $x\in A$. This
shows in particular that diffotopy is an equivalence
relation. \mn A locally convex algebra $A$ is called
contractible, if its identity map is diffotopic to 0, i.e. if
there is a continuous homomorphism $A\to A(0,1]$ such that
the composition with the evaluation map $A(0,1]\to A$ at 1 is
the identity map of $A$. The simplest examples of
contractible algebras are the algebras $\Cz(0,1]$ and the
algebra $t\Cz[t]$ of polynomials without constant term. The
required map $A\to A(0,1]$ maps in both cases the generator
$t$ to the function $s\to st$ for $s\in (0,1]$. \mn Below we
will also apply the notion of diffotopy to algebras $A$ over
$\Cz$ which we regard as a locally convex space with the
topology given by all possible semi-norms on $A$ (the fine
topology). In this case $\Cz[0,1]\hat {\otimes} A$ is simply
the algebraic tensor product $\Cz[0,1]\otimes A$.

\subsection{The algebra of smooth compact operators}
\label{compact} The $m$-algebra $\mathcal K$ of
``smooth compact operators" consists of all
$\Nz\times\Nz$-matrices $(a_{ij})$ with rapidly decreasing
matrix elements $a_{ij}\in\Cz,\,i,j=0,1,2\dots$. The topology
on $\mathcal K$ is given by the family of norms
$p_n,\,n=0,1,2\dots$, which are defined by
$$
p_n\bigl((a_{ij})\bigl)=\mathop\sum\limits_{i,j}\,
|1+i|^n|1+j|^n\, |a_{ij}|
$$
Thus, $\mathcal K$ is isomorphic to the projective tensor
product $s\hot s$, where $s$ denotes the space of rapidly
decreasing sequences $a=(a_i)_{i\in \Nz}$. The topology on
$s$ is determined by the seminorms $\alpha_n((a_i))=\sum
|1+i|^n\, |a_{i}|$. It is equipped with the natural bilinear
form $\langle a|b\rangle =\sum a_i b_i$. In this description,
the matrix product in $\mathcal K$ is given by $(a_1\otimes
b_1)(a_2\otimes b_2)= a_1\otimes \langle
b_1\mathop{|}a_2\rangle b_2$ using the canonical bilinear
form. This shows that the $p_n$ are submultiplicative. Since
$s$ is isomorphic, as a locally convex space, to the Schwartz
space $\mathcal S(\Rz)$ which in turn is isomorphic to
$\Cz(0,1)$ we see that $\mathcal K \cong \Cz(0,1)\hot
\Cz(0,1)$. It is important to note that these identifications
can be chosen to be compatible with the canonical bilinear
forms on $s$, $\mathcal S (\Rz)$ and $\Cz(0,1)$,
respectively, and therefore also with the multiplication on
$\mathcal K\cong s\hot s$. \blemma There are isomorphisms
$s\cong \mc S(\Rz)\cong \Cz (0,1)$ which respect the
canonical bilinear forms on these spaces (where the bilinear
form on $\mc S(\Rz )$ and $\Cz (0,1)$ is given by $\langle
f|g\rangle =\int fg$).\elemma \bproof The isomorphism between
$s$ and $\mc S (\Rz)$ is given by expansion of $f\in \mc
S(\Rz )$ with respect to the orthonormal basis of Hermite
polynomials, while an isomorphism between $\mc S(\Rz )$ and
$\Cz (-\pi/2,\pi/2)$ can for instance be obtained by
composition with the function arctan, cf. e.g. \cite{MV},
29.5.\eproof We call a homomorphism $\mathcal K\to \mathcal
K$ canonical if it is induced by a continuous linear map
$s\to s$ respecting the bilinear form. The bilinear form
preserving isomorphism $s\cong \Cz(0,1)$ and the induced
isomorphism $\mathcal K \cong \Cz(0,1)\hot \Cz(0,1)$ can be
used to prove the following lemma.
\begin{slemma}\label{homcomp} All canonical homomorphisms
$\mathcal K\to \mathcal K$ are diffotopic.
\end{slemma}
\bproof Let $\alpha$ be such a canonical homomorphism which
is induced by the bilinear form preserving map $\beta$. It is
clear that there exist continuous bilinear form preserving
linear maps $\beta_1,\beta_2 : \Cz(0,1)\to \Cz(0,1)$ such that
$\beta_1$ is diffotopic to the identity while $\beta_2$ is
diffotopic to $\beta$ (through bilinear form preserving maps)
and such that $\langle \beta_1(x)\mathop{|}\beta_2 (y)\rangle
=0$ for all $x,y$ in $\Cz(0,1)$. Put
$$\varphi_t (x) = \cos (t)\beta_1(x) + \sin (t)\beta_2(x)$$
Then $\varphi_t$ preserves the bilinear form for all $t$ and
induces a diffotopy between a canonical homomorphism
diffotopic to id and one diffotopic to $\alpha$.\eproof There is an obvious bilinear form preserving isomorphism $s\cong s\hot s$ which induces an isomorphism $\theta :\mc K\to \mc K\hot\mc K$. By the previous lemma, this isomorphism
is diffotopic to the inclusion map
$\iota:\mathcal{K}\to\mathcal{K}\hat{\otimes}\mathcal{K}$
that maps $x$ to $e_{00}\otimes x$ (where $e_{00}$ is the
matrix with elements $a^{ij}$, for which $a^{ij}=1$, whenever
$i=j=0$, and $a^{ij}=0$ otherwise).
\subsection{The smooth Toeplitz
algebra}\label{sToeplitz} The elements of the algebra
$\cC^\infty S^1$ can be written as Laurent series in the
generator $z$\, (defined by $z(t)=t$, $t\in S^1\subset\Cz)$.
The coefficients of these series are rapidly decreasing. Thus
$$
\cC^\infty(S^1)=\Big\{\,\mathop\sum\limits_{k\in\Zz} \,
a_kz^k\,\; \Big| \: \mathop\sum\limits_{k\in\Zz} |a_k |\,
|k|^n\, < \infty \ \mbox{for all $n\in \Nz$} \Big\}
$$
whence a natural isomorphism with the space $s=s(\Zz)$ of
rapidly decreasing sequences. Submultiplicative norms that
determine the topology are given by
$$
q_n\left(\sum a_kz^k\right)=\sum |1+k|^n\, |a_k|
$$
This topology on $\cC^\infty S^1$ is obviously the same as the
one used in \ref{df} for algebras of differentiable
functions.\\ As a topological vector space, the smooth
Toeplitz algebra $\mathcal T$ may be defined as the direct sum
$\mathcal{T}=\mathcal{K}\oplus \cC^\infty(S^1)$. To describe
the multiplication in $\mathcal{T}$, we let $\mc T$ act on
$s=s(\Nz)$. The elements of $\mc K$ act in the natural way
while the linear combinations $\sum a_kz^k$ act by truncated
convolution: For $a=\sum_{k\in\Zz} a_kz^k$ and $\xi
=(\xi_i)_{i\in\Nz}\in s$ we define $a\xi$ by $$(a\xi)_i =
\sum_{k+j=i}a_k\xi_j\;,\qquad i\in \Nz$$ In this way we
represent $\mc T$ as a subspace of End$(s)$ and define the
multiplication on $\mc T$ as the one inherited from End$(s)$.
We denote by $v$ and $v^*$ the images of $z$ and $z^{-1}$ in
$\mc T\subset \End (s)$. Then $v$ and $v^*$ are just the
right and left shift operators on the sequence space $s$. We
have $v^*v=1$ and $1-vv^*=e$ where e is the projection onto
the first basis vector of $s$.\\ If $p_n$ are the norms on
$\mathcal{K}$ defined in \ref{compact} and $q_n$ the norms on
$\cC^\infty(S^1)$ defined above, it is easy to see that each
norm $p_n \oplus q_n$ is submultiplicative on $\mathcal T
=\mathcal{K}\oplus \cC^\infty(S^1)$ with the multiplication
introduced above. Obviously, $\mathcal{K}$ is a closed ideal
in $\mathcal{T}$ and the quotient $\mathcal{T}/\mathcal{K}$ is
$\cC^\infty(S^1)$.
\begin{sproposition}\label{growth}
$\mathcal T$ is the universal unital $m$-algebra generated by
two elements $v$ and $v^*$ satisfying the relation $v^*v=1$
such that the assignment $z^k \mapsto v^k$ and $z^{-k}\mapsto
(v^*)^k$ for $k\in \Nz$, extends to a continuous linear map
$\mathcal C^\infty (S^1) \to \mathcal T$.
\end{sproposition}For the proof see \cite{CuDoc}, 6.1.\\
The smooth Toeplitz algebra (without this specific locally
convex topology) had been used before also by Karoubi,
\cite{Karo}).
\section{Extensions and classifying maps}
Let $V$ be a complete locally convex space. Consider the
algebraic tensor algebra
$$
T_{alg}V = V\,\oplus\, V\!\!\otimes\!\! V\,\oplus\,
V^{\otimes^{3}} \oplus\,\dots
$$
with the usual product given by concatenation of tensors.
There is a canonical linear map $\sigma:V\to T_{alg}V$ mapping
$V$ into the first direct summand. We equip $T_{alg}V$ with
the locally convex topology given by the family of all
seminorms of the form $\alpha\circ \varphi$, where $\varphi$
is any homomorphism from $T_{alg}V$ into a locally convex
algebra $B$ such that $\varphi\circ\sigma$ is continuous on
$V$, and $\alpha$ is a continuous seminorm on $B$.
Equivalently, we equip $T_{alg}V$ with the locally convex
topology determined by all submultiplicative sequences of
seminorms (see section \ref{lconv}) $(\alpha_1,\alpha_2,\ldots
)$ where all $\alpha_i$ are continuous on $V\subset
T_{alg}V$. We further denote by $TV$ the completion of
$T_{alg}V$ with
respect to this locally convex structure.\\
A more concrete description of the locally convex structure
on $TV$ has been determined by Valqui, \cite{Val}, who used
this free extension to generalize the proof of excision in
periodic cyclic homology from $m$-algebras to general locally convex algebras.
\begin{itemize}
\item $TV$ is a locally convex algebra.
\item Every continuous linear map $s: V \to B$ into a locally
convex algebra $B$ induces a continuous homomorphism
$\gamma_s : TV \to B$ defined by $$\gamma_s (x_1\otimes
\ldots \otimes x_n ) = s(x_1)s(x_2)\ldots s(x_n)$$ In
particular, $TV$ is contractible.
\item If $A$ is a locally convex algebra, then the canonical
homomorphism $TA\to A$ induced by the identity map of $A$ is
continuous. \end{itemize} For any locally convex algebra $ A$
we have the natural extension
\begin{equation}\label{UniExt} 0 \to J A \to T A \stackrel{\pi}{\to}  A \to 0.
\end{equation} Here
$\pi$ maps a tensor $x_1\otimes x_2\otimes\ldots \otimes x_n$
to $x_1x_2 \ldots x_n \in  A$ and $J A$ is defined as Ker
$\pi$. This extension is (uni)versal in the sense that, given
any extension $0 \to  I \to  E \to  B \to 0$ of a locally
convex algebra $B$, admitting a continuous linear splitting,
and any continuous homomorphism $\alpha :A\to B$, there is a
morphism of extensions \bgl\label{cla}
\begin{array}{ccccccccc} 0 &\to &  J A & \to & T A  & \to &   A
& \to & 0\\[0.1cm]
&      &  \quad\downarrow {\scriptstyle\gamma}& &
\quad\downarrow{\scriptstyle\tau} & & \quad\downarrow
{\scriptstyle\alpha} & &\\[0.1cm] 0
& \to &  I & \to &  E & \to &  B & \to & 0
\end{array}
\egl The map $\tau :T A\to  E$ is obtained by choosing a
continuous linear splitting $s: B\to  E$ in the given
extension and mapping $x_1\otimes x_2\otimes\ldots \otimes
x_n$ to $s'(x_1)s'(x_2) \ldots s'(x_n) \in E$, where
$s'\defeq s\circ \alpha$. Then $\gamma$ is the restriction of
$\tau$.\\ Choosing $0\to JB\to TB\to B\to 0$ in place of the
second extension $0\to I\to E\to B\to 0$, we see that $A \to
JA$ is a functor, i.e. any homomorphism $\alpha :A\to B$
induces a homomorphism $J(\alpha) :JA\to JB$.
\begin{definition} The map $\gamma
:J A \to  I$ in the commutative diagram (\ref{cla}) is called
the \emph{classifying map} (associated to the extension
together with the homomorphism $\alpha$). In the special case
where $B=A$ and $\alpha = \id$, we simply call $\gamma$ the
classifying map for the extension $0 \to I \to
 E \to  A \to 0$.\end{definition}\label{defclass}
If $s$ and $\bar{s}$ are two different continuous linear
splittings for a given extension, then so is $s_t =ts+(1-t)
\bar{s}$ for each $t$ in $[0,1]$. This differentiable family
of splittings induces a diffotopy between the classifying
maps defined from $s$ or $\bar{s}$. Therefore the classifying
map is unique up to homotopy.\\ Classifying maps can also be
defined for linearly split extensions of higher length of the
form
$$0\to I\to E_1\to E_2\to \ldots \to E_n\to A\to 0$$
We consider such an extension as a complex, denoting the
arrows (boundary maps) by $\pi_i$ and we say that it is
linearly split if there is a continuous linear map $s$ of
degree -1 such that $s\pi + \pi s = \id$. Every such
splitting $s$ induces a commutative diagram of the form
$$\begin{array}{ccccccccccccccc}0&\to
&I&\to& E_1&\to &E_2&\to& \ldots& \to &E_n&\to &A&\to& 0\\&
&\uparrow& & \uparrow& &\uparrow& & \ldots&  &\uparrow&
&\uparrow& &  \\0&\to &J^nA&\to&T(J^{n-1}A)&\to &T(J^{n-2}A)
&\to& \ldots& \to &TA&\to &A&\to& 0\end{array}$$ The leftmost
vertical arrow in this diagram is the \emph{classifying map}
for this $n$-step extension. It depends on $s$ only up to
diffotopy.
\section{Definition of the bivariant $K$-groups}\label{kkdef}
\begin{definition}\label{K-add} Let $ A$ and $ B$ be  locally convex  algebras.
For any continuous homomorphism $\varphi: A \to  B\,$,
we denote by $\langle \varphi \rangle$ the equivalence class
of $\varphi$ with respect to diffotopy and we set
$$
\langle  A , B\rangle =\{\langle \varphi \rangle |\,\varphi
\,\,\mbox{\rm is a continuous homomorphism }  A \to B\,\}
$$
\edefin Note that we work with diffotopy rather than with continuous homotopy since we want to construct a character into cyclic theory and since cyclic theory is invariant only under differentiable homotopies.\\ We denote by $\Cz(0,1)^n \cong \Cz(0,1)\hot \ldots
\hot\Cz(0,1)$ the algebra of infinitely differentiable
functions on the $n$-cube $[0,1]^n$ whose derivatives vanish on the
boundary. For every locally convex algebra $B$ and every
$n\geq 1$, there is a canonical map $\varphi :B(0,1)^n \oplus B(0,1)^n
\lori B(0,1)^n$, which is unique up to diffotopy. For any
second locally convex algebra $A$, this makes $\langle A ,\,
B(0,1)^n\, \rangle$ into an abelian group in the usual way by putting $\langle\alpha\rangle + \langle\beta\rangle= \varphi\circ (\alpha\oplus\beta)$.\\
Moreover, for each $k$ in $\Nz$, there is a canonical map $$\langle J^k A ,\, \mathcal K \hat{\otimes} B(0,1)^k\,
\rangle \to \langle J^{k+1} A ,\, \mathcal K \hat{\otimes}
B(0,1)^{k+1}\, \rangle$$ mapping the diffotopy class of
$\alpha$ to the diffotopy class of $\alpha '$, where
$\alpha'$ is defined by the commutative diagram
\[
\begin{array}{ccccccccc} 0 &\to &  J^{k+1}A & \to & TJ^k A  &
\to &   J^k A & \to & 0\\[0.1cm]
&      &  \qquad\downarrow {\scriptstyle \alpha'}& &
\quad\downarrow{\scriptstyle\tau }& & \quad\downarrow {\scriptstyle
\alpha} & &\\[0.1cm] 0 &
\to &  B'(0,1)^{k+1} & \to &  B'(0,1)^k[0,1) & \to &
B'(0,1)^k & \to & 0
\end{array}
\]
where $B' = \mc K\hot B$.
\begin{definition} Let $ A$ and $ B$ be locally convex algebras. We define
$$ kk^{\rm alg}( A ,\,  B\,) =
\lim_{\mathop{\lori}\limits_{k}} \langle J^k A ,\, \mathcal K
\hat{\otimes} B(0,1)^k\, \rangle
$$
\end{definition}
\bremark An alternative and basically equivalent definition of
$kk^{\rm alg}$ could be based on equivalence classes of
extensions of arbitrary length, rather than on their
classifying maps. Such an alternative approach has been
discussed, in the case of $C^*$-algebras, in
\cite{Thom}.\eremark We can also define abelian groups
$kk^{\rm alg}_n$ for each $n$ in $\Zz$ by
$$ kk^{\rm alg}_n( A ,\,  B\,) =
\lim_{\mathop{\lori}\limits_{k}} \langle J^{k-n} A ,\,
\mathcal K \hat{\otimes} B(0,1)^k\, \rangle
$$
where the inductive limit is taken over all $k\in \Nz$ such
that $k-n\geq 0$.\\
Thus, for $n$ in $\Nz$,
$$kk^{\rm alg}_n(A,B)=kk^{\rm alg}(J^nA,B)\, ,\qquad
kk^{\rm alg}_{-n}(A,B)=kk^{\rm alg}(A,B(0,1)^n)$$ There is a
natural product $kk^{\rm alg}(A_1,A_2)\times kk^{\rm
alg}(A_2,A_3)\lori kk^{\rm alg}(A_1,A_3)$ defined in the
following way. Assume that two elements $a$ in $kk^{\rm
alg}(A_1,A_2)$ and $b$ in $kk^{\rm alg}(A_2,A_3)$ are
represented by $\alpha :J^nA_1\to \mc K\hot A_2(0,1)^n$ and
$\beta :J^mA_2\to \mc K\hot A_3(0,1)^m$, respectively. We
define their product $a\cdot b$ in $kk^{\rm alg}(A_1,A_2)$ by
the following composition of maps
$$ J^m(J^nA_1)\quad \mathop{\lori}\limits^{J^m(\alpha)}\quad
J^m(\mc K\hot A_2(0,1)^n) \,\to\, \mc K\hot (J^m A_2)(0,1)^n
\quad \mathop{\lori} \limits^{\beta(0,1)^n}\quad \mc K\hot
A_3(0,1)^{n+m}$$ The arrow in the middle is obtained from the
natural map $J(D\hot E)\to J(D)\hot E$ (choosing $E=\mc
K(0,1)$) that exists for any two locally convex algebras $D$
and $E$, while the maps $J^m(\alpha)$ and $\beta(0,1)^n$ are
the natural maps induced by $\alpha$ and $\beta$,
respectively. The product is well defined (i.e. does not
depend on the choice of representatives for $a$ and $b$ in
the inductive limits defining $kk^{\rm alg}$) because of the
following lemma \blemma\label{well}(cf. \cite{Thom},C 3.10
and C 3.11). For any locally convex algebra $D$ we denote by
$\rho_D:JD\to D(0,1)$ the classifying map for the suspension
extension $0\to D(0,1)\to D[0,1)\to D\to 0$.\\ (a) Let
$$0\lori I\to E \mathop{\lori} \limits^{\pi} A\lori 0$$ be a
linearly split extension of locally convex algebras and
$\alpha$ its classifying map. Then the two maps $J^2A\to
I(0,1)$ defined by the compositions
$J^2A\mathop{\lori}\limits^{J(\alpha)}JI
\mathop{\lori}\limits^{\rho_I} I(0,1)$ and
$J^2A\mathop{\lori}\limits^{\rho_{JA}} (JA)(0,1)
\mathop{\lori}\limits^{\alpha(0,1)}I(0,1)$ are
diffotopic.\\
(b) The natural map $\rho_{JA}:\,J^2A\to (JA)(0,1)$ is
diffotopic to the composition $J^2A\to J(A(0,1))\to
(JA)(0,1)$. \elemma \bproof (a) The two compositions are
classifying maps for the two 2-step extensions in the first
and last row of the following diagram
\[
\begin{array}{ccccccccccc} 0 &\to &  I(0,1) & \to & I[0,1)  &
\to &   E&\to &A
& \to & 0\\[0.1cm]
&      &  \parallel& &
\downarrow & & \downarrow & &\downarrow& &\\[0.1cm] 0
& \to &  I(0,1) & \to &  E[0,1) & \to &  Z_\pi & \to & A&\to
& 0\\[0.1cm]&      &  \parallel& &
\uparrow & & \uparrow & &\uparrow& &\\[0.1cm]
0 & \to &  I(0,1) & \to &  E(0,1) & \to &  A[0,1) & \to &
A&\to & 0
\end{array}
\]
where $Z_\pi$ is the mapping cylinder. They are diffotopic
since they also are classifying maps for
the extension in the middle row.\\
(b) The first map is the classifying map for the first row in
the following commutative diagram \[
\begin{array}{ccccccccccc} 0 &\to &  (JA)(0,1) & \to & (TA)(0,1)
& \to &   A[0,1)&\to &A
& \to & 0\\[0.1cm]
&      &  \quad\uparrow{\scriptstyle\varphi} & &
\uparrow & & \uparrow& & \uparrow & &\\[0.1cm] 0
&\to &  J(A(0,1)) & \to & T(A(0,1)) & \to &   A[0,1)&\to &A &
\to & 0
\end{array}\]
while the second map is the composition of the classifying
map for the second row composed with the first vertical arrow
$\varphi$ in the diagram - thus also a classifying map for the
first row. \eproof Note that every continuous homomorphism
$\alpha :A\to B$ induces an element $kk(\alpha)$ in $kk^{\rm
alg}(A,B)$ and that $kk(\beta\circ\alpha)=kk(\alpha)\cdot
kk(\beta)$. We denote by $1_A$ the element $kk(\id_A)$
induced by the identity homomorphism of $A$. It is a unit in
the ring
$kk^{\rm alg}(A,A)$.\\
More generally, of course, we obtain a bilinear product
$kk_s^{\rm alg}(A_1,A_2)\times kk_t^{\rm alg}(A_2,A_3)\to
kk_{s+t}^{\rm alg}(A_1,A_3)$. Every $n$-step extension $$0\to
I\to E_1\to E_2\to \ldots \to E_n\to A\to 0$$ defines via its
classifying map an element in $kk_{-n}^{\rm alg}(A,I)$. Given
two such extensions of length $n$ and $m$, the product of the
corresponding elements in $kk_{-n}^{\rm alg}$ and
$kk_{-m}^{\rm alg}$ is represented by the classifying map for
the Yoneda product (concatenation) of the two extensions.
\begin{remark}The definition that we gave for $kk^{\rm alg}$
furnishes an abelian group also without the tensor product by
$\mc K$ in the second variable and the (bilinear) product can
be defined in exactly the same way. We refer to \cite{Thom}
for a thorough discussion of the theory obtained without
stabilization or only stabilization by finite matrices on the
right hand side, in the case of $C^*$-algebras. If we
stabilize however by $\mc K$ on the right hand side we have a
second possibility to define the addition in $kk^{\rm
alg}(A,B)$ (this second definition was used in \cite{CuDoc}).
For two continuous homomorphisms $\alpha , \beta : D \to
\mathcal K \hat{\otimes} E$ we define the direct sum
$\alpha\oplus\beta$\index{$\alpha\oplus\beta$} as
$$\alpha\oplus
\beta=\left(\begin{array}{cc}
\alpha&0\\0&\beta\end{array}\right) : D \lori M_2(\mathcal K
\hat{\otimes}E)\cong \mathcal K \hat{\otimes}E$$ With the
addition defined by $\langle\alpha\rangle +
\langle\beta\rangle = \langle\alpha\oplus\beta\rangle$ the
set $\langle D, \mathcal K \hat{\otimes}E\rangle$ of
diffotopy classes of homomorphisms from $D$ to $\mathcal K
\hat{\otimes}E$ is an abelian semigroup with $0$-element
$\langle 0 \rangle$. It is easy to see that this
definition of addition coincides with the one given above in
the case where $E=\mc K\hot B(0,1)^n$ for $n\geq 1$.
\end{remark}
In the case where we stabilize by $\mc K$ on the right hand
side, we can also use the following ``basic lemma'' from
\cite{CuDoc} to show that the product is well defined.
\begin{lemma}Assume given a commutative diagram of the form
$$
\begin{array}{ccccccccc}
 && 0 && 0 && 0 && \\
 && \downarrow && \downarrow && \downarrow && \\
 0 & \to &  I & \to &  A _{01} & \to &  A _{02} & \to & 0 \\
 && \downarrow && \downarrow && \downarrow && \\
 0 & \to & A _{10} & \to & A _{11} & \to & A _{12} & \to & 0 \\
 && \downarrow && \downarrow && \downarrow && \\
 0 & \to & A _{20} & \to & A _{21} & \to & B & \to & 0 \\
 && \downarrow && \downarrow && \downarrow && \\
 && 0 && 0 && 0 &&
\end{array}
$$ where all the rows and columns represent extensions of
locally convex algebras with continuous linear splittings.\\
Let $\gamma_+$ and $\gamma _-$ denote the classifying maps
$J^2 B \to I$ for the two extensions of length 2
$$ 0 \to I \to A _{10} \to A _{21} \to B \to 0$$
$$ 0 \to I \to A _{01} \to A _{12} \to B \to 0$$
associated with the two edges of the diagram. Then $\gamma _+
\oplus \gamma _- :\, J^2 B \to M_2(I)$ is diffotopic to 0.
\end{lemma}\begin{proof}
We assume first that the extension $0\to A_{20}\to A_{21}\to
B\to 0$ is trivial (i.e. admits a homomorphism splitting) and
show that then $\gamma_+\sim 0$  ($\gamma_- =0$ is obvious in
this case). By definition, $\gamma_+$ is constructed as a
classifying map from the
following diagram\\
$$
\begin{array}{ccccccccc}
0 & \rightarrow & I & \rightarrow & A_{01} & \rightarrow &
A_{02} & \rightarrow & 0\\
 & & \uparrow & & \uparrow & & \uparrow \alpha & &\\
0 & \rightarrow & J^2B & \rightarrow & TJB & \rightarrow & JB
& \rightarrow & 0
\end{array}
$$
where $\alpha$ is the classifying map for the extension
$$
0 \rightarrow A_{02} \rightarrow A_{12} \rightarrow B
\rightarrow  0
$$
The commutative diagram
$$
\begin{array}{ccccccccc}
0 & \rightarrow & A_{01} & \rightarrow & A_{11} & \rightarrow
& A_{21} & \rightarrow & 0\\
 & & \;\downarrow {\scriptstyle\pi } & & \downarrow & &
 \downarrow & &\\
0 & \rightarrow & A_{02} & \rightarrow & A_{12} & \rightarrow
& B & \rightarrow & 0
\end{array}
$$
with the splitting homomorphism $\beta :\, B\to A_{21}$ shows
that $\alpha: JB \rightarrow A_{02}$ can be lifted to a
homomorphism $\alpha' : JB \rightarrow A_{01}$ such that
$\pi\circ\alpha' =\alpha$. In fact, choose $\alpha'$ as the
classifying map for
$$
\begin{array}{ccccccccc}
0 & \rightarrow & A_{01} & \rightarrow & A_{11} & \rightarrow
& A_{21} & \rightarrow & 0\\
 & & & & & & \;\uparrow {\scriptstyle\beta } & &\\
 & & & & & & B & &
\end{array}
$$
It follows that $\gamma_+ \sim 0$.\mn Returning now to the
general case, consider the following commutative diagram
$$
\begin{array}{ccccccccc}
  && 0 && 0 && 0 && \\[2mm]
  && \downarrow && \downarrow && \downarrow && \\[2mm]
  0 & \to &  I & \to & I & \to &  0 & \to & 0 \\[2mm]
  && \downarrow && \downarrow && \downarrow && \\[2mm]
  0 & \to & A_{01} + A _{10} & \to & A _{11} & \to & B &
  \to & 0 \\[2mm]
  && \downarrow && \downarrow && \downarrow && \\[2mm]
  0 & \to & A_{02} + A _{20} & \to & A _{11} / I & \to &
  B & \to & 0 \\[2mm]
  && \downarrow && \downarrow && \downarrow && \\[2mm]
  && 0 && 0 && 0 &&
\end{array}
$$
From the first step it follows that the classifying map for
the (linearly split) extension
$$
0 \to I \to A_{01} + A_{10} \to A_{11} / I \to B \to 0
$$
is diffotopic to $0$.\\
An obvious rotation argument shows that the composition of
the classifying map $J^2B\to I$ for the extension
$$
\begin{array}{ccccccccc}
0 & \rightarrow & I & \rightarrow & A_{01} + A_{10} &
\rightarrow & A_{02} \oplus A_{20} & \rightarrow & 0\\
 & & & & & & \uparrow  & &\\
 & & & & & & JB & &
\end{array}
$$
with the inclusion $I \hookrightarrow M_2(I)$ is diffotopic to
the map $x \longmapsto \left(\begin{array}{ll}
\gamma_+(x) & 0\\
0 & \gamma_-(x)
\end{array}\right)$
where $\gamma_+ ,\, \gamma_- :\,J^2B\to I$ are the classifying
maps associated with the diagrams
$$
\begin{array}{ccccccccc}
0 & \rightarrow & I & \rightarrow & A_{01} & \rightarrow &
A_{02} & \rightarrow & 0\\
 & & & & & & \quad \uparrow {\scriptstyle \eta_+ }& &\\
 & & & & & & JB & &
\end{array}
$$
and
$$
\begin{array}{ccccccccc}
0 & \rightarrow & I & \rightarrow & A_{10} & \rightarrow &
A_{20} & \rightarrow & 0\\
 & & & & & & \quad \uparrow {\scriptstyle \eta_- } & &\\
 & & & & & & JB & &
\end{array}
$$
and $\eta_+,\, \eta_-$ are classifying maps for the
extensions $0 \to A_{02} \to A_{11} \to B \to 0$ and $0 \to
A_{20} \to A_{11} \to B \to 0$.
\end{proof}
This lemma is much stronger than Lemma \ref{well} in the case
where the right hand side is stabilized by $\mc K$. It also
shows that the product of the elements of $kk_{-1}^{\rm alg}$
induced by the classifying maps for the extensions $ 0 \to I
\to A _{01} \to A _{02}\to 0$ and $0\to A _{02}\to A _{12}
\to B \to 0$ is the negative of the product of the elements
of $kk_{-1}^{\rm alg}$ induced by the classifying maps for the
extensions $ 0 \to I \to A _{10} \to A _{20}\to 0$ and $0\to
A _{20}\to A _{21} \to B \to 0$. \\ $kk^{\rm alg}_*$ is best
viewed as an additive category whose objects are locally
convex algebras and whose (graded) morphism sets between
objects $A$ and $B$ are $kk^{\rm alg}_*(A,B)$. The map
$\alpha\to kk(\alpha)$ is a functor from the category of
locally convex algebras to the category $kk^{\rm alg}_*$.
\bprop Let $A$ be a locally convex algebra. There is an
invertible element in $kk^{\rm alg}(JA,A(0,1))$ giving rise to
an isomorphism $kk^{\rm alg}(JA,B)\cong kk^{\rm
alg}(A(0,1),B)$. \eprop \bproof The natural map $J(A(0,1))\to
(JA)(0,1)$ can be viewed as defining an element in $kk^{\rm
alg}(A(0,1),JA)$. The classifying map for the suspension
extension
$$0\lori A(0,1)\lori A[0,1)\lori A\lori 0$$ can be viewed as
defining an element in $kk^{\rm alg}(JA,A(0,1))$. The
products of these two elements give the identity elements in
$kk^{\rm alg}(JA,JA)$ and $kk^{\rm alg}(A(0,1),A(0,1))$\eproof
As a consequence we see that $$kk^{\rm alg}(A(0,1),B(0,1))\cong kk^{\rm alg}(JA,B(0,1))\cong kk^{\rm
alg}(A,B)$$
\section{Long exact sequences}
Let $\alpha :A\to B$ be a continuous homomorphism between
locally convex algebras. The mapping cone $C_\alpha \subset
A\oplus B[0,1)$ is defined to be the set of pairs
$(x,f)\in A\oplus B[0,1)$ such that $f(0)=\alpha (x)$.
\begin{lemma}\label{puppe1}
Let $D$ be a locally convex algebra and $\alpha: A\to B$ a
continuous homomorphism between locally convex algebras
\begin{enumerate}
\item[(a)] The sequence
$$
kk^{\rm alg}_*(D,C_\alpha)\mathop{\lori}\limits^{\cdot
kk(\pi)}kk^{\rm alg}_*( D, A) \mathop{\lori}\limits^{\cdot
kk(\alpha)}kk^{\rm alg}_*( D, B)
$$
is exact. Here $\pi:C_\alpha\to A$ denotes the projection
onto the first summand and $\cdot kk(\pi)$ multiplication on
the right by $kk(\pi)$.
\item[(b)] The sequence in (a) can be extended to an exact
sequence
$$
\begin{array}{llllll}
\mathop{\lori}\limits^{\cdot kk(\pi (0,1))}&kk^{\rm alg}_*(
D, A (0,1))& \mathop{\lori}\limits^{\cdot
kk(\alpha(0,1))}&kk^{\rm alg}_*( D, B\,(0,1))&\to
& \\[6pt] & kk^{\rm alg}_*( D, C_\alpha)&
\mathop{\lori}\limits^{\cdot kk(\pi)}&kk^{\rm alg}_*( D, A)
\mathop{\lori}\limits^{\cdot kk(\alpha)}&kk^{\rm alg}_*( D,
B)&
\end{array}
$$
\end{enumerate}
\end{lemma}
\bproof (a) follows directly from the definition of the
mapping cone. In fact, let $\alpha': X\to Y$ be any continuous
homomorphism between locally convex algebras and $\varphi$ a
homomorphism from $D$ to $X$. Then $\alpha\circ\varphi$ is
diffotopic to 0 if and only if there is a homomorphism
$\tilde{\varphi}:D\to C_{\alpha'}$ such that $\varphi =
\pi\circ\tilde{\varphi}$.\\
(b) follows by a standard argument from the iteration of the
construction under (a). One uses the fact that the mapping
cone $C_\pi$ for the projection $\pi: C_\alpha\to A$ is
diffotopic to $ B (0,1)$ and the following commutative diagram
$$
\begin{array}{ccc}
C_\pi  & \mathop{\lori}\limits^{\pi'} & C_\alpha\\
\uparrow     &  & \parallel  \\
 B\,(0,1) & \mathop{\lori}\limits^{\iota} & C_\alpha
\end{array}
$$
In this diagram $\iota$ is the inclusion of $ B(0,1)$ into
the second component of $C_\alpha$ and the first vertical
arrow is a diffotopy equivalence (it maps $f\in  B\, (0,1)$
to $(\iota f, 0) \in C_\pi \subset C_\alpha \oplus  A
[0,1))$. \mn Similarly, the mapping cone $C_\iota$ for
$\iota: B(0,1)\to C_\alpha$ is contained in $ A(0,1)\oplus
 B([0,1)\times [0,1))$. The projection $C_\iota \to
 A (0,1)$ is a diffotopy equivalence making the
following diagram commute
$$
\begin{array}{ccc}
 C_\iota  & \lori &  B(0,1)\\
\downarrow     &  & \parallel  \\
 A\, \,(0,1) & \mathop{\lori}\limits^{\alpha(0,1)} &
 B\, (0,1)
\end{array}
$$\eproof
\begin{lemma}\label{puppe2}  Let $\alpha: A\, \to  B\, $ and $ D$
be as in \ref{puppe1}
\begin{enumerate}
\item[(a)] The sequence
$$
kk^{\rm alg}_\ast(C_\alpha,  D)
\mathop{\longleftarrow}\limits^{kk(\pi) \cdot} kk^{\rm
alg}_\ast ( A, D) \mathop{\longleftarrow}\limits^{kk(\alpha)
\cdot} kk^{\rm alg}_\ast( B,
 D)
$$
is exact.
\item[(b)] The sequence in (a) can be extended to a long exact
sequence of the form
$$
\begin{array}{llllll}
\mathop{\longleftarrow}\limits^{ kk(\pi (0,1))\cdot}&kk^{\rm
alg}_*( A (0,1), D)&
\mathop{\longleftarrow}\limits^{kk(\alpha(0,1))\cdot}&kk^{\rm
alg}_*(
B(0,1),  D)&\leftarrow & \\[6pt] & kk^{\rm alg}_*(C_\alpha , D)&
\mathop{\longleftarrow}\limits^{kk(\pi)\cdot}&kk^{\rm alg}_*(
A, D)
\mathop{\longleftarrow}\limits^{kk(\alpha)\cdot}&kk^{\rm
alg}_*( B , D)&
\end{array}
$$
\end{enumerate}
\end{lemma}
\bproof Let $\varphi:J^n A\to \mathcal K\hat{\otimes}
D(0,1)^n$ represent an element $z$ of $kk^{\rm alg}(A,D)$
such that $\varphi\circ J^n(\pi)$ is diffotopic to 0. To
simplify notation, we write $D'$ for $\mathcal K\hat{\otimes}
D(0,1)^n$. The fact that $\varphi\circ J^n(\pi)$ is
diffotopic to 0 means that there is a homomorphism $\gamma$
such that the following diagram is commutative
$$
\begin{array}{ccc}
 D'[0,1) &\lori  & D'\\[3pt]
 \uparrow {\scriptstyle\gamma} & &\uparrow {\scriptstyle\varphi} \\[4pt]
J^{n} C_\alpha & \mathop{\lori}\limits^{J^n(\pi)} & J^{n} A
\end{array}
$$
The restriction of $\gamma$ defines a map $\gamma'$ such that
the following diagram is commutative
$$
\begin{array}{ccccccc}
0\lori &  D'(0,1) & \lori &  D'
[0,1) &\lori  & D' &  \lori  0\\[3pt]
& \uparrow {\scriptstyle\gamma'} && \uparrow {\scriptstyle\gamma}
 & &\uparrow {\scriptstyle\varphi} \\[4pt]
0\lori &  J^n(B(0,1)) & \lori & J^{n} C_\alpha &
\mathop{\lori}\limits^{J^n(\pi)} & J^{n} A\,   & \lori 0
\end{array}
$$
Let $\kappa_B :J^{n+1}B\to J^n(B(0,1))$ denote the natural map
and let $z'$ denote the element of $kk^{\rm alg}(B,D)$
defined by
$\gamma'\circ\kappa_B$.\\
Define $\psi$ by the commutative diagram
$$
\begin{array}{cccc}
& & &J^n(A(0,1))\qquad \\[4pt]
&  &  &\downarrow {\scriptstyle\mu_t} \\[4pt]
 C_{J^n(\pi)}  & \subset && J^nA (-1,0]\oplus J^n(C_\alpha)
 \\[5pt]
\downarrow {\scriptstyle\psi}  &  & &\qquad \downarrow
{\scriptstyle \varphi(-1,0]
\oplus \gamma} \\[5pt]
 D' (-1,1)\qquad & \subset &&
 D'\,(-1,0]\oplus D' \,[0,1)
\end{array}
$$
There is a differentiable family of homomorphisms $\mu_t :
J^n(A(0,1))\lori C_{J^n(\pi)}$ such that $\mu_0$ is the
canonical map $J^n(A(0,1))\to (J^n A)(0,1)$ composed with the
natural inclusion $$(J^n A)(0,1)\cong (J^n A)(-1,0)\lori
C_{J^n(\pi)}$$ while $\mu_1$ is the composition of
$J^n(\alpha (0,1))$ with the map $J^n(B(0,1))\to J^n
(C_\alpha)$
induced by the inclusion $B(0,1)\to C_\alpha$ (we use here
a straightforward modification of the standard homotopies for
maps $A(0,1)\to C_\pi$ which also appear in the proof of
\ref{puppe1} (b).\\
Let now finally $\kappa_A :J^{n+1}A\to J^n(A(0,1))$ denote the
natural map. Then $\psi\circ\mu_t\circ \kappa_A$ defines a
diffotopy between the image of $\varphi$ in $\langle
J^{n+1}A,\mathcal K\hat{\otimes} D(0,1)^{n+1}\rangle$ and the
composition $\gamma'\circ \kappa_B\circ J^n(\alpha)$.\\
Thus $z = kk(\alpha)\cdot z'$ as desired. \eproof
\begin{proposition}\label{Cq}
Let $0\to  I\to  A\, \mathop{\lori}\limits^q  B\, \to 0$ be a
linearly split extension and $e: I\to C_q$ the inclusion map
defined by $e:x\mapsto (x,0)\in C_q \subset A\,  \oplus B\,
[0,1)$. Then $kk(e)$ is an invertible element in $kk^{\rm
alg}_0( I, C_q)$.
\end{proposition}
\bproof We show that the inverse is given by the element $u$
of $kk^{\rm alg}_0(I, C_q)$ determined by the classifying map
$JC_q \to I$ for the extension
$$
0\lori I(0,1)\lori A\, [0,1) \lori C_q \lori 0
$$
The commutative diagram
 $$
 \begin{array}{ccccccc}
 0\lori & I(0,1) & \lori &  A[0,1) & \lori & C_q & \lori 0\\
 & \uparrow & & \uparrow & & \uparrow {\scriptstyle e} &\\
 0\lori &  I(0,1) & \lori &  I[0,1) & \lori &  I & \lori 0
 \end{array}
 $$
 with the obvious maps shows that $kk(e)\cdot u=1_I$.\\
 On the other hand, let $\Delta :A[0,1)\to A[0,1)[0,1)$ be the
 diagonal map defined by $(\Delta f)(x,y)= f(\abs{(x,y)})$ for
 $(x,y)\in [0,1]\times [0,1]$. Let also
 $\psi :A[0,1)[0,1)\to C_q[0,1)$ be the map induced by the
 obvious quotient map $A[0,1)\to C_q$.\\
 Defining $e'$ as $\psi\circ \Delta$ we obtain a commutative
 diagram
$$
 \begin{array}{ccccccc}
 0\lori & C_q(0,1) & \lori & C_q[0,1) & \lori & C_q & \lori 0
 \\[2pt]
 & \qquad\uparrow {\scriptstyle e(0,1)}& & \uparrow {\scriptstyle
 e'}& & \parallel &\\[3pt]
 0\lori &  I(0,1) & \lori &  A[0,1) & \lori & C_q & \lori 0
 \end{array}
$$ which shows that $u\cdot kk(e)=1_{C_q}$.\eproof
\begin{theorem}
Let $ D$ be any locally convex algebra. Every extension
admitting a continuous linear section
$$
E:\, 0\to  I\stackrel{i}{\lori}  A\, \mathop{\lori}\limits^q
 B\,\to 0
$$
induces exact sequences in $kk^{\rm alg}( D,\cdot\,)$ and
$kk^{\rm alg}(\,\cdot\,,  D)$ of the following form: \bglnoz
\ldots\lori kk^{\rm alg}_{1}( D,A) \stackrel{\,\cdot
kk(q)}{\lori} kk^{\rm alg}_{1}(
D, B)\qquad\qquad\\[4pt]\lori kk^{\rm alg}_0( D,I)\stackrel{\,
\cdot kk(i)} {\lori} kk^{\rm alg}_0( D, A) \stackrel{\,\cdot
kk(q)}{\lori} kk^{\rm alg}_0( D, B)\lori \ldots \eglnoz and
\bglnoz \ldots \lole kk^{\rm alg}_{-1}( A,D) \stackrel{\,\cdot
kk(q)}{\lole} kk^{\rm alg}_{-1}(
B,D)\qquad\qquad\\[4pt]\lole kk^{\rm alg}_0( I,D)\stackrel{\,
\cdot kk(i)} {\lole} kk^{\rm alg}_0(A,D) \stackrel{\,\cdot
kk(q)}{\lole} kk^{\rm alg}_0(B,D)\lole \ldots \eglnoz The
given extension $E$ defines a classifying map $J B \to  I$
and thus an element of $kk^{\rm alg}_{-1}( B, I)$, which we
denote by $kk(E)$. The boundary maps are (up to a sign) given
by right and left multiplication, respectively, by this class
$kk(E)$.
\end{theorem}
\bproof The long exact sequences follow from a combination of
\ref{puppe1}, \ref{puppe2} and \ref{Cq}. To determine the
boundary maps, consider the mapping cylinder $Z_q \subset
A\oplus B[0,1]$ consisting of all pairs $(x,f)\in A\oplus
B[0,1]$ such that $f(0)=q(x)$.\\
The commutative diagram
$$
 \begin{array}{ccccccc}
 0\lori & B(0,1) & \lori & B[0,1) & \lori & B & \lori 0
 \\[2pt]
 & \downarrow& & \downarrow & & \parallel &\\[3pt]
 0\lori &  C_q & \lori &  Z_q & \lori & B & \lori 0\\[2pt]
 & \quad\uparrow {\scriptstyle e}& & \uparrow & & \parallel &\\[2pt]
 0\lori &  I & \lori &  A & \lori & B & \lori 0
 \end{array}
$$ shows that $kk(E)\cdot kk(e)$ equals $kk(j)$ in $kk^{\rm alg}_0(B(0,1),C_q)\cong
kk^{\rm alg}_1(B,C_q)$ where $j: B(0,1)\to C_q$ is the natural
inclusion. By construction, the boundary map in \ref{puppe1},
\ref{puppe2} is multiplication by $kk(j)$. On the other hand,
the identification between $I$ and $C_q$ in \ref{Cq} is given
by multiplication with $kk(e)$. \eproof A locally convex
algebra $C$ is called $k$-contractible if it is isomorphic in
the category $kk^{\rm alg}$ to 0, or equivalently if $kk^{\rm
alg}_0(C,C)=0$ (which implies that $kk^{\rm
alg}_n(A,C)=kk^{\rm alg}_n(C,A)=0$ for all $n$ and $A$). Of
course, every contactible algebra is $k$-contractible.
\bcor\label{kcontr} Let $E:\;0\to I\to C\to B\to 0$ be a
linearly split extension where $C$ is $k$-contractible. Then
the element $kk(E)\in kk^{\rm alg}_{-1}(B,I)$ is invertible,
i.e. there is $u\in kk^{\rm alg}_{-1}(I,B)$ such that
$kk(E)\cdot u =1_B$ and $u\cdot kk(E) =1_I$.\\ More
generally, if $\gamma$ is the classifying map for a linearly
split $n$-step extension
$$0\to I\to C_1\to \ldots \to C_n\to B\to 0$$
where $C_1,\ldots,C_n$ are $k$-contractible, then the element
$kk(\gamma)\in kk^{\rm alg}_0(J^n(B),I)\cong kk^{\rm
alg}_{-n}(B,I)$ is invertible. \ecor \bproof Since every
$n$-step extension is a Yoneda product of 1-step extensions
and $kk(\gamma)$ is the product of the $kk^{\rm alg}$-elements
corresponding to the corresponding classifying maps, it
suffices to prove the assertion for the first case of
one-step extensions.\\ In this case the product by $kk(E)$
gives by \ref{puppe1} the boundary map $kk^{\rm
alg}_{1}(D,B)\to kk^{\rm alg}_0 (D,I)$ associated to the
extension $E$ for all $D$, thus also for $D=I$. Since, $C$ is
$k$-contractible, this boundary map must be an isomorphism.
Therefore there is an element $u\in kk^{\rm alg}_{1}(I,B)$
such that $u\cdot kk(E)=1_I$.\\ Similarly, the boundary map
$kk^{\rm alg}_{1}(I,B)\to kk^{\rm alg}_0(B,B)$, given by left
multiplication by $kk(E)$ in \ref{puppe2} is an isomorphism.
Therefore there is an element $u'\in kk^{\rm alg}_{1}(I,B)$
such that $kk(E)\cdot u' =1_B$. \eproof
\section{Quasihomomorphisms}
Let $\alpha$ and $\bar{\alpha}$ be two homomorphisms $A\to E$
between locally convex algebras. Assume that $B$ is a closed
subalgebra of $E$ such that $\alpha (x)-\bar{\alpha}(x)\in B$
and $\alpha(x)B\subset B$, $B\alpha(x)\subset B$ for all
$x\in A$. We call such a pair $(\alpha,\bar{\alpha})$, a
quasihomomorphism from $A$ to $B$. We construct an element
$kk(\alpha,\bar{\alpha})$ of $kk^{\rm alg}(A,B)$ in the
following way. Define $\alpha', \bar{\alpha}':A\to A\oplus E$
by $\alpha'(x)=(x,\alpha(x)),
\,\bar{\alpha}'=(x,\bar{\alpha}(x))$ and $E'$ as the subalgebra
of $E\oplus A$ generated by all elements $\alpha'(x),\
x\in A$ and by $0\oplus B$. Let $$D = \left\{f\in
E'[0,1]\mathop{\big|} \exists x\in A, f(0)=\alpha'(x),
f(1)=\bar{\alpha}'(x)\right\}$$ Then the classifying map for
the linearly split extension
$$0\to B(0,1)\to D\to A\to 0$$
defines the element $kk(\alpha,\bar{\alpha})$ of $kk^{\rm alg}(A,B)$.\\
This construction is obviously invariant under diffotopy: If
$\alpha$ and $\bar{\alpha}$ are two homomorphisms $A\to
E[0,1]$ such that $\alpha (x)-\bar{\alpha}(x)\in B[0,1]$ and
$\alpha(x)B[0,1]\subset B[0,1]$, $B[0,1]\alpha(x)\subset
B[0,1]$ for all $x\in A$ then $kk(\alpha_t,\bar{\alpha_t})$ is
constant. \bremark\label{quasi} It suffices to assume that
$\alpha (x)-\bar{\alpha}(x)\in B$ and $\alpha(x)B\subset B$,
$B\alpha(x)\subset B$ for all $x$ in a generating subset $G$
of $A$.\eremark \blemma\label{dif} Let $\varphi =\alpha -
\bar{\alpha}$. If
$\varphi(x)\bar{\alpha}(y)=\bar{\alpha}(y)\varphi(x)=0$ for
all $x,y,\in A$, then $\varphi$ is a homomorphism and
$kk(\alpha,\bar{\alpha})=kk(\varphi)$.\elemma \bproof The
first assertion is obvious and the second one follows from the
following commutative diagram
$$
 \begin{array}{ccccccc}
 0\lori & B(0,1) & \lori & D & \lori & A & \lori 0
 \\[2pt]
 & \qquad\uparrow {\scriptstyle \varphi(0,1)}& & \uparrow
 {\scriptstyle \psi}& & \parallel &\\[3pt]
 0\lori &  A(0,1) & \lori &  A[0,1) & \lori & A & \lori 0
 \end{array}$$
where $\psi(f)(t) = \bar{\alpha}(f(0))+\varphi(f(t))$. \eproof
\bprop Let $E$ be a locally convex algebra containing $\mc K\hot A$ as a closed ideal and $\alpha,\bar{\alpha} :A\to E[0,1]$ a pair of
homomorphisms such that $\alpha (x)
- \bar{\alpha} (x)\in\mc K\hot A[0,1)$ for all $x$ in $A$.\\
Assume further that $\varphi =\alpha_0 -\bar{\alpha}_0$ is a
homomorphism, satisfies the assumptions of \ref{dif} and is
diffotopic to the natural inclusion map $\iota :A\to \mc
K\hot A$. Then $A$ is $k$-contractible in the sense of
\ref{kcontr}.\eprop \bproof We have $0= kk(\alpha_1
,\bar{\alpha}_1)=kk(\alpha_0 ,\bar{\alpha}_0) = kk(\iota)$.
On the other hand $kk(\iota)$ is an invertible element in
$kk^{\rm alg}(A,\mc K\hot A)$. \eproof
\section{Morita invariance}
We discuss in this section a rather weak form of Morita
equivalence which is however sufficient for our purposes in
later sections. \bdefin\label{mor} Let $A$ and $B$ be locally
convex algebras. A Morita context from $A$ to $B$ is a locally
convex algebra $E$ containing $A$ and $B$ as subalgebras and
containing two sequences $(\xi_i)$ and $(\eta_j)$ in $E$
satisfying
\begin{itemize}
\item $\eta_j A\xi_i\subset B$ for all $i,j$.
\item the family $(\eta_ja\xi_i)$ is rapidly decreasing for each
$a\in A$
(i.e. $\alpha(\eta_ja\xi_i)$ is rapidly decreasing in $i$
and $j$ for any continuous seminorm $\alpha$ on B).
\item $(\sum \xi_i\eta_i )a =a$ for all $a\in A$
(i.e the partial sums converge in $A$) \end{itemize}\edefin
Given a Morita context, we can
define a continuous homomorphism $A\to \mc K\hot B$ by
$$ a \mapsto  \big(\eta_i a\xi_j\big)$$ where we write the
elements of $\mc K\hot B$ as $\Nz\times\Nz$-matrices. Thus,
in particular, we obtain an element of $kk^{\rm alg}_0(A,B)$
which we denote by $kk((\xi_i),(\eta_j))$.
\blemma\label{Morita} Every Morita context ($(\xi_i)$,
$(\eta_j)$) from $A$ to $B$ induces an element
$kk((\xi_i),(\eta_j))$ of $kk^{\rm alg}(A,B)$. If
$((\xi'_l),(\eta'_k))$ is a Morita context from $B$ to $A$
realized in the same locally convex algebra $E$ and if in
addition $A\xi_i\xi'_l\subset A$, $\eta'_k\eta_jA\subset A$,
for all $i,j,k,l$, then
$$kk((\xi_i),(\eta_j))\cdot kk((\xi'_l),(\eta'_k))= 1_A$$
If also $B\xi'_l\xi_i\subset B$, $\eta_k\eta'_j B\subset B$
for all $i,j,k,l$, then $A$ and $B$ are $kk^{\rm
alg}$-equivalent.\elemma \bproof The composition of the two
induced maps $$A\to \mc K\hot B\to \mc K\hot\mc K\hot A\cong
\mc K\hot A$$ maps $a\in A$ to the $\Nz^2\times\Nz^2$-matrix
$$\Big(\eta'_k\eta_j a \xi_i\xi'_l\Big)$$ This map is
diffotopic to the canonical inclusion of $A$ into $\mc K\hot
A$ via the homotopy that maps $a\in A$ to the $L\times
L$-matrix $(\hat{\eta}_\alpha(t) a \hat{\xi}_\beta(t))$ where
$L=\Nz^2\cup \{0\}$ and
$$\begin{array}{ccc} \hat{\xi}_0(t)=\cos t\,1&&
\hat{\xi}_{il}(t)=\sin t\;\xi_i\xi'_l\\[4mm]
\hat{\eta}_0(t)=\cos t\,1&& \hat{\eta}_{kj}(t)=\sin
t\;\eta'_k\eta_j
\end{array}$$ The conditions $A\xi_i\xi'_l\subset A$,
$\eta'_k\eta_jA\subset A$ ensure that the elements
$(\hat{\eta}_0(t) a \hat{\xi}_{il}(t))$ and
$(\hat{\eta}_{kj}(t) a \hat{\xi}_0(t))$ are in $A$.\eproof
\bremark\label{mmor} A typical example of two algebras that
are Morita equivalent in the sense of \ref{Morita} are
algebras $A$ and $B$ where $B$ is a subalgebra of $\mc K\hot
A$ such that $(e_{00}\otimes 1)\,B\,(e_{00}\otimes 1) =
e_{00}\otimes A$ and such that $B(e_{j0}\otimes 1)\subset B$
and $(e_{0j}\otimes 1)B\subset B$ for all $j$. In fact, we
can choose $\xi_0=\eta_0= e_{00}\otimes 1$ (the other
$\xi_i,\eta_j$ all being $0$) and $\xi'_l = e_{l0}\otimes 1,
\, \eta'_k = e_{0k}\otimes 1$ in this case. \eremark A
situation closely related to Morita invariance arises in the
following way. Let $V$ be a complete locally convex space
equipped with a continuous non-zero bilinear form
$\langle\;|\;\rangle$. Then $\mc M=V\hot V$ becomes a locally
convex algebra. Assume now that $\psi :V\to s$ is a
continuous linear map that preserves the bilinear form. It
induces a continuous homomorphism $\alpha_\psi :\mc M \to \mc
K$. We have the following \blemma The map
$\alpha_\psi\otimes\id :\mc M\hot A\to \mc K\hot A$ induces a
$kk^{\rm alg}$-equivalence for each locally convex algebra
$A$.\\ Its inverse is given by the natural inclusion $A\to
\mc M\hot A$ mapping $x\in A$ to $e\otimes x$ where $e$ is a
one-dimensional projection in $\mc M$. \elemma \bproof Let
$\xi\in s$ and $\eta \in V$ be vectors such that $\langle\xi
|\xi\rangle=\langle\eta |\eta\rangle=1$. Consider the family
of maps $\psi_t : V\to V\hot s \oplus V\hot s$ defined by $$
\psi_t (x) = (\cos t (x\otimes \xi), \sin t (\eta\otimes \psi
(x)))$$ for $t\in [0,\pi /2]$. It induces a differentiable
family of continuous homomorphisms $\alpha_{\psi_t}:\mc M \to
M_2(\mc M\hot \mc K)$ which is for $t=0$ the natural
inclusion of $\mc M$ into $\mc M\hot \mc K$ (in the upper
left corner of $M_2$) and for $t=\pi/2$, $\alpha_\psi$
composed with the natural inclusion of $\mc K$ into $\mc
M\hot \mc K$ (in the lower right corner of $M_2$). \eproof
\section{Bott periodicity}
Recall that we defined the smooth
Toeplitz algebra $\mc T$ as the universal unital $m$-algebra
generated by two elements $v$ and $v^*$ satisfying the
relation $v^* v=1$ such that the canonical linear map $\mc
C^\infty S^1 \to \mc T$ becomes continuous. It fits into an
extension
$$0 \to \mc K \to \mc T \to \mc C^\infty S^1\to 0$$
We denote by $e$ the idempotent $1-vv^*$ in $\mc T$.
\begin{lemma}\label{pair} (compare \cite{CuK}, 4.2)
There are unique continuous homomorphisms $\varphi, \varphi ':
\mc T\to \mc T\hat{\otimes} \mc T$, such that
$$
\begin{array}{l}
\varphi(v)= v(1-e)\otimes1 + e\otimes v \qquad
\varphi(v^*)=(1-e)v^*\otimes1 + e\otimes v^*
\\[2pt]
\varphi '(v)=v(1-e)\otimes1 + e\otimes 1 \qquad
\varphi'(v^*)=(1-e)v^*\otimes1 + e\otimes 1
\end{array}
$$
These two homomorphisms are diffotopic through a diffotopy
$\psi_t: \mc T\to \mc T\hat{\otimes} \mc T, \, t\in
[0,\pi/2]$, for which $\psi_t(x)-\varphi(x)\in \mc
K\hat{\otimes} \mc T$ holds for all $t\in[0,\pi/2],x\in  \mc
T$.
\end{lemma}
\bproof We show, that $\varphi$ and $\varphi'$ are both
diffotopic to $\psi$, where
$$
\psi(v)= v\otimes 1 \qquad \psi(v^*)=v^*\otimes 1
$$
We write $e_{ij}$ for the matrix units in $\mc K\subset \mc T$. With this notation
we write linear combinations $e_{ij}\otimes x, \: 0\le i,j \le
n-1,\: x \in \mc T$ as $n\times n$-matrices with entries in $
\mc T$. For $t\in [0,\pi/2]$ we
put
$$
\begin{array}{l}
u_t=v^2v^{*2}\otimes 1 + \left(\begin{array}{cc}e_{00}+\cos
t\:vv^* & \sin t\: v \\ -\sin t\: v^* & \cos t\,1
\end{array}\right)
\\[10pt]
u'_t=v^2v^{*2}\otimes 1\:+ \:\left(\begin{array}{cc} \cos
t\,1& \sin t \,1\\ -\sin t \,1& \cos t\,1
\end{array}\right)
\end{array}
$$
Then $u_t$ and $u'_t$ are obviously invertible elements in $
\mc T\hat{\otimes} \mc T$ and lie in fact in the subalgebra
generated by $ \mc K\hat{\otimes} \mc T$ and $ 1\hat{\otimes}
\mc T$. By the characterization of $\mc T$ as a universal
algebra in \ref{growth} there exist continuous homomorphisms
$\varphi_t,\varphi'_t: \mc T\to \mc T\hat{\otimes} \mc T$,
such that
$$
\begin{array}{l}
\varphi_t(v)= u_t\,(v\otimes1) \qquad
\varphi_t(v^*)=(v^*\otimes1)\,u_{t}^{-1}
\\[2pt]
\varphi'_t(v)= u'_t\,(v\otimes1) \qquad
\varphi'_t(v^*)=(v^*\otimes1)\,u_{t}'^{-1}
\end{array}
$$
For the verification of the necessary growth conditions in
\ref{growth} we refer to \cite{CuDoc}, 6.2. \eproof We may
consider the algebra $\Cz (0,1)$ as the subalgebra of $\mc
C^\infty S^1$ of functions vanishing, together with all
derivatives at the point 1. Denote by $\mc T_0$ the preimage
of $\Cz (0,1)$ in the extension $$0 \to \mc K \to \mc T \to
\mc C^\infty S^1\to 0$$ so that we obtain an extension
$$0 \to \mc K \to \mc T_0 \to \Cz (0,1)\to 0$$ Consider also
the algebra $\mc T'_0 \supset \mc T_0$ defined as the kernel
of the natural map $\mc T\to \Cz$ mapping $v$ and $v^*$ to 1.
The commutative diagram $$
 \begin{array}{ccccccc}
 0\lori & \mc K & \lori &  \mc T'_0 & \lori & D
  & \lori 0\\
 & \uparrow & & \uparrow & & \uparrow &\\
 0\lori &  \mc K & \lori &  \mc T_0 & \lori &  \Cz (0,1)
 & \lori 0
 \end{array}
$$ where $D=\left\{f\in \mc C^\infty
 (S^1)\,\big|\,f(1)=0\right\}$ and the last vertical arrow is a
diffotopy equivalence shows that $kk^{\rm alg}(\mc T_0,\mc
T_0)\cong kk^{\rm alg}(\mc T_0,\mc T'_0)$. Moreover, the split
extension $$0\to \mc T'_0\to \mc T\to\Cz\to 0$$ shows that
$kk^{\rm alg}(\mc T_0,\mc T'_0)\cong\quad kk^{\rm alg}(\mc
T_0,\mc T)\oplus kk^{\rm alg}(\mc T_0,\Cz)$. Tensoring
everything with an arbitrary locally convex algebra $A$ shows
that the natural inclusion map $\mc T_0\hot A \to \mc T\hot
A$ induces an injective map $kk^{\rm alg}_0(\mc T_0\hot A,\mc
T_0\hot A)\to kk^{\rm alg}_0(\mc T_0 \hot A,\mc T\hot A)$ for
each locally convex algebra $A$. \bprop\label{T0} For each
locally convex algebra $A$, the algebra $\mc T_0 \hot A$ is
$k$-contractible (i.e. isomorphic to 0 in the category
$kk^{\rm alg}$).\eprop \bproof It suffices to show that the
natural inclusion map $kk^{\rm alg}_0(\mc T_0\hot A,\mc
T_0\hot A)\to kk^{\rm alg}_0(\mc T_0 \hot A,\mc T\hot A)$ is
0. Denote by $\bar{\varphi}$ the continuous homomorphism $\mc
T \to \mc T \hot \mc T$ defined by
$$\bar{\varphi}(v)= v(1-e)\otimes1\qquad
\bar{\varphi}(v^*)=(1-e)v^*\otimes1$$ Using the diffotopy
$\psi_t$ from the previous lemma, we consider the restriction
of the pair $(\psi_t,\bar{\varphi})$ to homomorphisms from
$\mc T_0$ to $\mc T\hot \mc T$. It induces for all $t$ the
same element $kk(\psi_t,\bar{\varphi})$ in $kk^{\rm alg}(\mc
T_0, \mc K\hot \mc T)\cong kk^{\rm alg}(\mc T_0, \mc T)$.
However, for $t=0$ (where $\psi_0=\varphi$) we have $\psi_0
-\bar{\varphi} = \iota$ where $\iota$ is the natural
inclusion map $\mc T_0\to \mc K\hot\mc T$, while for $t=1$
(where $\psi_1=\varphi'$) we have $\psi_{1} -\bar{\varphi} =
0$ on $\mc T_0$.\\The same argument applies if we tensor $\mc
T_0$ by $A$ and all homomorphisms by $\id _A$. \eproof
\btheo\label{Bott} (Bott periodicity) The extension
$$E:\; 0\to \mc K\hot A \to \mc T_0\hot A\to A(0,1)\to 0
$$ defines an invertible element $\varepsilon$ in
$kk^{\rm alg}_{-1}(A(0,1),A)\cong kk^{\rm alg}_{-2}(A,A)$. As
a consequence $kk^{\rm alg}_n (A,B)\cong kk^{\rm
alg}_{n-2}(A,B)$ for all $A,B$ and $n$. \etheo\bproof This
follows from \ref{kcontr} and \ref{T0} .\eproof
\begin{remark}\label{gleich} The fact that the inverse
Bott map $\varepsilon :J^2A\to \mc K\hot A$ associated as a
classifying map to the extension$$0\to \mc K\hot A\to\mc
T_0\hot A \to\Cz [0,1)\hot A\to A\to 0$$ defines an
invertible element, shows
that the definition of $kk^{\rm alg}$ that we use in this
article gives the same result as the analogue of the formally
different definition given in \cite{CuDoc}:
$$kk^{\rm alg}_\ast(A ,\, B\,)
=\lim_{\mathop{\lori}\limits_{k}} \langle J^{2k+\ast}
 A ,\, \mathcal K \hat{\otimes} B\, \rangle
$$ where the inductive limit is taken using the
$\varepsilon$-map from $J^{k+2}A$ to $\mc K\hot J^kA$. In
fact, there are obvious natural maps from the theory defined
that way into $kk^{\rm alg}$ defined as above and vice versa,
which are inverse to each other.\\ Note also that conceivable
different versions of the Toeplitz algebra would not work for
the proof of Bott periodicity above. A natural candidate for
instance would be the enveloping $m$-algebra of the unital
algebra generated algebraically by two elements $v$ and $v^*$
satisfying the relation $v^*v=1$. This algebra would not work
since the quotient by the natural ideal $\mc K'$ generated by
$e=1-vv^*$, which is a different (smaller than $\mc K$)
completion of the algebra $M_\infty$ of finite matrices of
arbitrary size, is too small to admit a homomorphism from
$\Cz (0,1)$ into it. This is closely related to the fact,
studied by Thom in \cite{Thom}, that stabilization in the
definition $kk^{\rm alg}$ on the right hand side by algebras
smaller than $\mc K$ leads to theories which are not
periodic, such as connective $K$-theory. \end{remark}
\begin{theorem}\label{sixterm}
Let $ D$ be any locally convex algebra. Every extension
admitting a continuous linear section
$$
E:\, 0\to  I\stackrel{i}{\lori}  A\, \mathop{\lori}\limits^q
 B\,\to 0
$$
induces exact sequences in $kk^{\rm alg}( D,\cdot\,)$ and
$kk^{\rm alg}(\,\cdot\,,  D)$ of the following form: \bgl
\label{1exact}
\begin{array}{ccccc}
kk^{\rm alg}_0( D,  I) & \stackrel{\,\cdot kk(i)}{\lori} &
kk^{\rm alg}_0( D,
 A) &\stackrel{\,\cdot kk(q)}{\lori} & kk^{\rm alg}_0(
D, B)
\\[3pt]
\uparrow  & & & & \downarrow
\\[2pt]
kk^{\rm alg}_1( D,  B) & \stackrel{\,\cdot
kk(q)}{\longleftarrow} & kk^{\rm alg}_1( D,  A) &
\stackrel{\,\cdot kk(i)}{\longleftarrow} & kk^{\rm alg}_1( D,
I)
\end{array}
\egl and
\bgl \label{2exact}
\begin{array}{ccccc}
kk^{\rm alg}_0( I, D) & \stackrel{kk(i)\cdot
\,}{\longleftarrow} & kk^{\rm alg}_0(  A,  D) &
\stackrel{kk(q)\cdot \,}{\longleftarrow} & kk^{\rm alg}_0( B,
D)
\\[3pt]
\downarrow  & & & & \uparrow
\\[2pt]
kk^{\rm alg}_1( B,  D) & \stackrel{kk(q)\cdot\,}{\lori} &
kk^{\rm alg}_1( A,
 D) &\stackrel{kk(i)\cdot\,}{\lori} & kk^{\rm alg}_1(
I,  D)
\end{array}
\egl The given extension $E$ defines a classifying map $J B
\to  I$ and thus an element of $kk^{\rm alg}_1( B, I)$, which
we denote by $kk(E)$. The vertical arrows in (\ref{1exact})
and (\ref{2exact}) are (up to a sign) given by right and left
multiplication, respectively, by this class $kk(E)$.
\end{theorem}
\section{The bivariant Chern-Connes character}
With every pair of locally convex algebras $A$, $B$ we can
associate the bivariant cyclic homology groups $HP_*(A,B)$,
$\star = 0,1$ which are in fact complex vextor spaces. The
bivariant theory $HP_*$ has properties analogous to those of
$kk^{\rm alg}_*$. In particular it is diffotopy invariant in
both variables, invariant under tensoring one of the
arguments by $\mc K$ and, most importantly, it has six-term
exact sequences exactly like those in \ref{sixterm}, replacing
there $kk^{\rm alg}$ by $HP$, \cite{CQInv}, \cite{Val},
\cite{Mey}, \cite{CV} where the boundary maps are given by an
element in $HP_1$ determined by the given
extension.\\Therefore the standard extensions
$$0\to JA\to TA\to A\to 0$$
$$0\to A(0,1)\to A[0,1)\to A\to 0$$ determine isomorphisms
$HP_1(A,JA)\cong HP_0(A,A)$, $HP_1(A(0,1),A)\cong HP_0(A,A)$
and invertible elements $\alpha_A\in HP_1(A,JA)$ and
$\beta_A\in HP_1(A(0,1),A)$. By iteration we have canonical
invertible elements $\alpha_A^n\in HP_n(A,J^nA)$ and
$\beta_A^n\in HP_n(A(0,1)^n,A)$. Let moreover $\gamma^n$ be
the element of $HP_0(J^nA,A(0,1)^n)$ induced by the canonical
homomorphism $J^nA\to A(0,1)^n$. The following relations are
then easily checked:
$$\alpha_A^n\cdot\gamma^n\cdot\beta_A^n =1_A\qquad
\gamma^n\cdot\beta_A^n\cdot\alpha_A^n =1_{J^nA}\qquad
\beta_A^n\cdot\alpha_A^n\cdot\gamma^n =1_{A(0,1)}$$ where
$1_B$ denotes the identity element of $HP_0(B,B)$ (induced by
the identity map of $B$). \btheo\label{Chern} There are
natural maps $ch:kk^{\rm alg}_s(A,B)\to HP_s(A,B)$ such that
$$ch(\varphi\cdot\psi)=ch(\varphi)\cdot ch(\psi)$$ in
$HP_{s+t}(A,C)$ for $\varphi\in kk^{\rm alg}_s(A,B)$ and
$\psi\in kk^{\rm alg}_t(B,C)$ and such that $ch(1_A)=1_A$ for
each $A$. Here $s$ and $t$ are in $\Zz$ and $HP_s$ is
considered to be 2-periodic in $s$.\etheo\bproof Ignoring
first the stabilization by $\mc K$ assume that $\varphi\in
kk^{\rm alg}_s(A,B)$ is represented by a homomorphism $\eta
:J^{n+s}A\to B(0,1)^n$. We then define
$$ch(\varphi)= \alpha_A^{n+s}HP_0(\eta)\beta_B^n$$
The relation $\alpha_A^n\cdot\gamma^n\cdot\beta_A^n =1_A$
mentioned above generalizes to
$\alpha_A^n\cdot\theta^n\cdot\beta_B^n =\theta$ if $\theta$
is an element of $HP_0(A,B)$ determined by a homomorphism
$A\to B$ and $\theta^n$ is the element of
$HP_0(J^nA,B(0,1)^n)$ determined by the homomorphism $J^nA\to
B(0,1)^n$ induced by $\theta$. Since the inductive limit in
the definition of $kk^{\rm alg}_s$ is taken exactly with
respect to the passage from homomorphisms $A\to B$ to
homomorphisms $J^nA\to B(0,1)^n$, this shows that $ch$ is
well defined. The compatibility with the product follows from
Lemma \ref{well}.\\Finally, the natural inclusion $B\to \mc
K\hot B$ induces an invertible element $\iota$ in $HP_0(B,\mc
K\hot B)$. The formula for $ch$ then becomes $$ch(\varphi)=
\alpha_A^{n+s}HP_0(\eta)\beta_B^n\,\iota^{-1}$$ if $\varphi$
is represented by the map $\eta :J^{n+s}A\to\mc K\hot
B(0,1)^n$. The stabilization by $\mc K$ does not lead to any
essential change in the arguments for well-definedness and
multiplicativity. \eproof It is important to note that the
Bott periodicity element $\varepsilon\in kk^{\rm alg}_2(A,A)$
is not mapped to $1$ under $ch$, but rather to $2\pi i$, cf.
\cite{CuDoc}, 6.10.
\section{The coefficient ring and the exterior product}
\bprop There is a bilinear graded commutative exterior product
$$kk^{\rm alg}_i(A_1,B_1)\times kk^{\rm alg}_j(A_2,B_2)\lori
kk^{\rm alg}_{i+j}(A_1\hot A_2,B_1\hot B_2)$$ denoted by
$(\alpha,\beta)\mapsto \alpha\times\beta$. \eprop \bproof
There are natural maps $J(A_1\hot A_2)\to JA_1\hot A_2$ and
$J(A_1\hot A_2)\to A_1\hot J A_2$. For all $s,t\in \Nz$ there
is a natural map $J^{s+t}(A_1\hot A_2)\to J^sA_1\hot J^t
A_2$.\eproof The exterior product $\times$ and the interior
product $\cdot$ of section \ref{kkdef} are compatible in the following way:
\blemma\label{exin} Let $\alpha_1\in kk^{\rm alg}_0(A_1,B_1)$,
$\alpha_2\in kk^{\rm alg}_0(A_2,B_2)$, $\beta_1\in kk^{\rm
alg}_0(B_1,C_1)$, $\beta_2\in kk^{\rm alg}_0(B_2,C_2)$. Then
$$(\alpha_1\times\alpha_2)\cdot
(\beta_1\times\beta_2)=(\alpha_1\cdot\beta_1)\times
(\alpha_2\cdot\beta_2)$$\elemma On the coefficient ring
$R\defeq kk^{\rm alg}_0(\Cz,\Cz)$ the exterior and interior
products obviously coincide. With this product $R$ is a
commutative unital ring. $kk^{\rm alg}_i(A,B)$ is a module
over $R$ for all $A,B$ and $i$ and by Lemma \ref{exin} the
composition product
$$kk^{\rm alg}_i(A,B)\times kk^{\rm alg}_j(B,C)\lori
kk^{\rm alg}_{i+j}(A,C)$$ is $R$-linear.\\
By Theorem \ref{Chern}, there is a unital ring homomorphism
$ch: R\to \Cz =HP_0(\Cz,\Cz)$. In particular $R$ is not zero!
\section{The Weyl algebra and its fundamental extension}
\label{fund} The Weyl algebra $W$ is the unital algebra
generated by two elements $x,y$ satisfying the relation
$[x,y]=1$. A fundamental reference for the properties of this
algebra is \cite{Dix}. It is easy to see that $W$ can not
carry any submultiplicative seminorm. In fact, the completion
of $W$ with respect to such a seminorm would be a unital
Banach algebra containing two elements $x,y$ such that
$xy-yx=1$. This relation shows that for the (non-empty!)
spectra we have Sp($xy$)=Sp($yx$)+1. On the other hand we
would have $\Sp(xy)\cup\{0\}=\Sp(yx)\cup\{0\}$. This shows
that the spectrum of $xy$ would have to be unbounded which is
impossible in a Banach algebra.\\ As explained above however,
$W$ is a locally convex algebra with the fine topology and we
will from now on always assume that $W$ carries this topology.
\mn To construct the fundamental extension, let $W'$ denote
the algebra with two generators $x',y'$ satisfying the
relations
$$
(x'y' - y'x') y' = y' \qquad x'(x'y' - y'x') = x'
$$
We equip $W'$ first with the fine locally convex topology
(later we will also consider other locally convex structures
on $W'$).\mn Let $f= x'y'- y'x'-1$ in the algebra
$\widetilde{W'}$ with unit adjoined. This element will play
an important role in the following computations. The relations
$fy' =x'f= 0$ imply that
$$fx'y' = x'y' f = f + f^2$$
\begin{lemma}\label{xy} We have
$$x'^m(x'y') = (m+x'y')x'^m,\quad [x', x'y'] = x'$$
for each $m\in \Nz$.
\end{lemma}
\begin{proof}
We have $x'(x'y') = x'(y'x'+1) = (1+x'y')x'$ whence the
second equation. The first identity follows from this by the
derivation rule.
\end{proof}
Let $\pi :W'\to W$ be the natural map and $I$ the kernel of
$\pi$. Then $I$ is the ideal generated in $W'$ by all
elements of the form $fa$ and $af$, for $a\in W'$. \mn We
consider natural representations of $W$ and $W'$ as operators
on the spaces of (arbitrary) functions with compact support
$V_{\Rz} = \mathcal F_c([0,\infty))$ and $V_{\Nz} = c_0(\Nz)$
where $\Nz=\{0,1,2,\ldots\}$ (the support condition is in fact
inessential). We define representations $\rho':W'\to \End
V_{\Rz}$ and $\rho:W\to \End V_{\Nz}$ by
$$\rho'(x') = \sqrt{t}\,U_{-1} \; \; \; \rho'(y') = U_1
\sqrt{t}$$ and
$$\rho(x) = \sqrt{N}\,U_{-1} \; \; \; \rho(y) = U_1 \sqrt{N}$$
where $U_1$ and $U_{-1}$ are right and left translation:
$$(U_1\xi)(s) = \xi(s-1)\qquad (U_{-1}\xi)(s) = \xi(s+1)$$ and
$\sqrt{t},\,\sqrt{N}$ are the square roots of the ``number
operators" in $\End V_{\Rz}$ and $\End V_{\Nz}$ defined by
$$(\sqrt{t}\xi)(s) = \sqrt{s}\xi(s)\qquad (\sqrt{N}\xi)(n) =
\sqrt{n}\xi(n)$$ for $\xi \in V_{\Rz}, V_{\Nz}$.\\ It is clear
that $\rho'(x')$, $\rho'(y')$, and $\rho(x)$, $\rho(y)$
satisfy the relations of $W'$ and $W$, respectively. In fact
$\rho$ is nothing but the standard representation of the
Heisenberg commutation relation by creation and annihilation
operators used in quantum mechanics.
\begin{lemma}\label{w}
Every element $w$ in $I$ has a unique representation as a
finite sum of the form
$$
w = \sum_{k,l} w_{kl} \; \; \; k,l = 0,1,2,\ldots
$$
with $w_{kl} = y'^{k} P_{kl}(x'y')fx'^{l}$ where $f$ is
defined as before \ref{xy}, $P_{kl}$ is a polynomial in
$\Cz[t]$ in general and $P_{00}$ is in $t\Cz [t]$.\\
Moreover, the representation $\rho '$ is faithful on $I$
\end{lemma}
\begin{proof}
Every element of $I$ is a linear combination of elements of
the form $af,\,fb$ or $afb$ with $a,b$ in $W'$. We know that
$fy' =x'f =0$ and, from Lemma \ref{xy},
$$
fx'^n y'^k = f((n-1) + x'y') x'^{n-1}y'^{k-1}
$$
$$ = \left\{\begin{array}{lll}
((n-1) + x'y')(n-2+x'y') & \ldots (n-k+x'y')fx'^{n-k} &n \ge k
\\[3mm]
0 & n < k &
\end{array} \right. $$
and similarly
$$
x'^n y'^k f = y'^{(k-n)}P(x'y')f
$$
for a certain polynomial $P$. This shows that $w$ can be
written as claimed.\\ To show that $\rho'$ is faithful on $I$
note that the elements $w_{kl}$ are mapped under $\rho'$ to
$z_{kl}=U_1^k Q_{kl}P_{kl}(t)\rho'(f)U_{-1}^l$ where $Q_{kl}$
is a non-zero product of functions of the form $\sqrt
{t+n},\,n\in\Nz$, that depends only on $k$ and $l$ (in fact
it is not difficult to write it down explicitly) and
$\rho'(f)$ is multiplication by the function $\hat{f}$ given
by
$$ \hat{f}(t) = \left\{\begin{array}{ll}
1-t &\quad t \leq 1
\\[2mm]
0 & \quad t\geq 1
\end{array} \right. $$
It is clear that the operator $\sum z_{kl}$ in $\End (V_\Rz)$
is 0 only if all $P_{kl}$ are 0.\\
This shows at the same time the injectivity of $\rho'$ and
the uniqueness of the $P_{kl}$.
\end{proof}
\blemma\label{I}The ideal $I$ is Morita equivalent in the
sense of \ref{Morita} to the algebra $t(1-t)\Cz [t]$ of
polynomials vanishing in 0 and 1. As a consequence $I$ is
$k$-equivalent to $\Cz (0,1)$.\elemma \bproof According to
Lemma \ref{w} there is an inclusion map $t(1-t)\Cz [t]\to I$
mapping $(1-t)P$, with $P\in t\Cz [t]$, to $fP(x'y')$.\\ To
construct the Morita context from $I$ to $t(1-t)\Cz [t]$,
think of $W'$ as represented in End $V_\Rz$. Using the
notation of \ref{mor} we choose $\xi_i=\xi_i '=U_1^i\chi$ and
$\eta_j=\eta'_j=\chi U_{-1}^j$ where $\chi$ is the
multiplication operator
by the characteristic function for the interval $(0,1]$.\\
From the extension $0\to (1-t)t\Cz[t]\to t\Cz[t]\to \Cz\to 0$
where $t\Cz[t]$ is contractible (see subsection \ref{df}), it
follows that $(1-t)t\Cz [t]$ is isomorphic in the category
$kk^{\rm alg}$ to $\Cz (0,1)$. \eproof
\begin{lemma}\label{faith}
The representation $\rho'$ is faithful and we have a
commutative diagram
$$\begin{array}{lll}
W' & {\mathop{\rightarrow}\limits^{\rho'}} & {\rm End}\; V_{\Rz}\\
\,\downarrow {\scriptstyle \pi} & &\quad \downarrow
{\scriptstyle\eta}
\\
W & {\mathop{\rightarrow}\limits^{\rho}} & {\rm End}\; V_{\Nz}
\end{array}$$
where $\eta$ is restriction to the invariant subspace of
$V_\Rz$ consisting of functions with finite support in $\Nz$ (which is
isomorphic to $V_{\Nz}$).
\end{lemma}
\begin{proof} The commutativity of the diagram is obvious.
Since $W$ is simple it shows that
$\Ker \rho' \subset I = \Ker \pi$. But $\rho '$ is faithful
on $I$ by \ref{w}\\
\end{proof}

\section{Triviality of $W'$ and
the $kk^{\rm alg}$-invariants for $W$ .}\label{W'}
Using the notation of \ref{pair} we define invertible elements
$u_t$ in $\mc T$ by
$$u_t = v^2 v^{*2}\: + \:\left(\begin{array}{cc} \cos t&
\sin t\\ -\sin t & \cos t
\end{array}\right)$$ where the $2\times 2$-matrix represents a
linear combination of the matrix units $e_{ij}$, $i,j=0,1$ in
$\mc K\subset \mc T$. We can define a differentiable family
$\varphi_t, t\in [0,\pi/2]$ of homomorphisms $W'\to
W'\hot\mathcal T$ by
$$\varphi_t(x')=x'\otimes v^*u^{-1}_t\qquad \varphi_t(y')=
y'\otimes
u_tv$$ and another homomorphism $\bar{\varphi}$ defined by
$$\bar{\varphi}(x')=x'\otimes (1-e)v^*\qquad \bar{\varphi}
(y')=y'\otimes v(1-e)$$ For each $t\in [0,\pi/2]$ and each
$z\in W'$ we have that
$$\varphi_t(z)-\bar{\varphi}(z)\quad\in\quad W'\hot\mathcal
K\subset W'\hot\mathcal T$$ We therefore obtain induced
elements $kk(\varphi_t,\bar{\varphi})$ in $kk^{\rm
alg}(W',W'\hot \mathcal K)$ which do, by homotopy invariance,
not depend on $t$. \blemma The linear space $B$ generated by
all elements of the form $$y'^n cx'^m\otimes e_{nm}\;,\quad
n,m = 0,1,2\ldots$$ where $c$ is a polynomial in $x'y'$ in
general, and a polynomial without constant term for $n=m=0$,
is a subalgebra of $W'\hot\mathcal K$.\\ For each $z\in W'$
we have
$$\varphi_0(z) -\bar{\varphi}(z) \in B\qquad
\varphi_0(z)B\subset B\qquad B\varphi_0(z)\subset B$$ \elemma
\bproof The first assertion follows from the fact that
$x'^ny'^n$ is a polynomial without constant term in $x'y'$
according to Lemma \ref{xy}.\\According to Remark \ref{quasi}
it suffices to check the second assertion for $z=x'$ and
$z=y'$ where it follows again immediately from Lemma
\ref{xy}.\eproof \blemma Let $f_0=1$ and $f_n \in \Cz[t],\,
n\geq 1$ denote the polynomial
$$ f_n (t) = t(t+1)(t+2)\ldots(t+n-1)$$
and let $\mathcal K _0\subset \mathcal K$ denote the algebra
of finite matrices.\\ The map $y'^n cx'^m\otimes e_{nm}\mapsto
cf_m\otimes e_{nm}$ defines an injective homomorphism $B\to
\Cz[t]\hot\mathcal K_0$. The image of $B$ under this injection
consists of all linear combinations of elements $cf_m\otimes
e_{nm}$ with $c\in \Cz[t]$ in general, and $c\in t\Cz[t]$ for
$n=m=0$. \elemma \bproof The identity $x'^ny'^n = f_n(x'y')$
from Lemma \ref{xy}, shows that the map in question is a
homomorphism. Injectivity follows from \ref{faith}. The rest is obvious. \eproof \btheo The algebra
$W'$ is $k$-contractible, i.e. $kk^{\rm alg}(W',W')=0$.\etheo
\bproof We show first that $B$ is $k$-contractible. Denote by
$B_0$ the kernel of the homomorphism $B\to \mathcal K_0$
induced by the evaluation $t\mapsto 1$ from $\Cz[t]\hot
\mathcal K_0$ to $\mathcal K_0$. The image $N\subset \mathcal
K_0$ of $B$ under this map consists of all linear
combinations of elements $e_{n0},\, n\geq 1$ and is therefore
nilpotent of order 2 ($N^2=0$) and thus also contractible. We obtain an
extension
$$0\to B_0\to B\to N\to 0$$ $B_0$ is obviously Morita
equivalent to $t\Cz [t]$ in the sense of \ref{Morita} (see
Remark \ref{mmor}) and therefore also
$k$-contractible. It follows that $B$ is $k$-contractible.\\
As a consequence, with the notation above,
$kk(\varphi_{\pi/2},\bar{\varphi})=kk(\varphi_0,\bar{\varphi})=0$
since these morphisms are products of the corresponding
elements in $kk^{\rm alg}(W',B)$ by the element $kk(j)$
induced by the inclusion $B\to W'\hot \mathcal K$ and since
$B$ is $k$-contractible. On the other hand $\varphi_{\pi/2}-
\bar{\varphi}$ satisfies the conditions of \ref{dif} and is a
homomorphism equal to the canonical embedding $W'\to
W'\hot\mathcal K$. It follows that $kk^{\rm
alg}(W',W')=0$.\eproof Note that exactly the same argument
shows that $W'$ is isomorphic to 0 in $HP_*$, i.e.
$HP_*(W',W')=0$.\\ We are now in a position to compute the
$kk^{\rm alg}$-invariants for the Weyl algebra $W$. \btheo
The Weyl algebra $W$ is isomorphic to $\Cz$ in $kk^{\rm
alg}$.\etheo \bproof Since $W'$ is $k$-contractible, the
element $kk(E)$ in $kk^{\rm alg}_1(W,I)$ associated with the
extension $0\to I\to W'\to W\to 0$ is invertible. By\ref{I},
$$kk^{\rm alg}_1(W,I)\cong kk^{\rm alg}_1(W,\Cz (0,1))
\cong kk^{\rm alg}_0(W,\Cz)$$ \eproof It follows also from this
result that $W$ is isomorphic, in the category $HP_*$ defined
by bivariant periodic cyclic homology, to $\Cz$ (invertible
elements are mapped by the character $ch$ to invertible
elements). The Hochschild and cyclic homology of $W$
however have been known for a long time, cf. \cite{FT}.\\
From the result above we can also deduce the $kk^{\rm alg}$-
invariants for the higher dimensional versions $W_n$ of $W$.
By definition $W_n$ is generated by elements $x_1,\ldots ,
x_n$ and $y_1,\ldots , y_n$ such that $[x_i,y_j]=\delta_{ij}1$
and $[x_i,x_j]=0$, $[y_i,y_j]=0$ for all $i,j$. Thus $W_n
=W\hot \ldots \hot W =W\otimes \ldots \otimes W$.\\ The
$n$-fold exterior product of the invertible element $kk(E)$ in
$kk^{\rm alg}_1(W,I)\cong$\linebreak $kk^{\rm alg}(W,\Cz)$
constructed above, by itself, gives an invertible element in
$kk^{\rm alg}(W_n,\Cz)$. It is represented by the $n$-step
extension
$$ 0 \to I^{\otimes n} \to I^{\otimes n-1}\otimes W'\to
\ldots \to W'\otimes W^{\otimes n-1}\to W^{\otimes n}\to 0$$
Using similar extensions one may also compute the $kk^{\rm
alg}$-invariants for further algebras of similar type.
Consider for instance the algebra $\mc D$ of differential
operators generated by the Fr\'{e}chet-algebra $\mc S (\Rz)$
of Schwartz functions together with the differentiation
operator $D$, satisfying the relation $[D,f]=f'$ for $f\in
\mc S (\Rz)$. We make it into a locally convex algebra by
defining a seminorm $\alpha$ on $\mc D$ to be continuous, if
the map $\mc S (\Rz)\ni f \mapsto \alpha(fD^n)$ is continuous
for each $n$. Using a similar extension as for $W$ we can
show that $kk^{\rm alg}_0(\Cz,\mc D)=0$, $kk^{\rm
alg}_1(\Cz,\mc D)=R$.
\section{Seminorms on $W'$}
In this section, we construct a different locally convex topology
on $W'$ in such a way that the completion $\bar{I}$ of the ideal
$I$ becomes an $m$-algebra (in fact $I$ will be Morita equivalent
to $\mc C^\infty_0 (0,1)$ with the canonical $m$-algebra
structure). This gives an extension $$0\to \bar{I}\to \bar{W}'\to
W\to 0$$ where the ideal is an $m$-algebra. This extension will
be used in section \ref{app} to obtain more precise information
on the $K$-theoretic invariants of $W$.\mn Recall from Lemma
\ref{xy} the formulas
$$x'^m(x'y') = (m+x'y')x'^m,\quad (x'y')y'^m= y'^m(m+x'y')$$
Putting as above $f= (x'y'-y'x')-1$ in $\widetilde{W}'$, we
obtain
$$
x'^k y'^lf= x'^{k-1} y'^{l-1}((l-1) + x'y')f =
x'^{k-1}y'^{(l-1)}(l+f)f
$$
By induction (on $l$) we obtain
$$x'^k y'^l f= \left\{\begin{array}{lll}
y'^{l-k}(l+f)(l-1+f) \ldots & (l+1-k+f)f&
l \ge k\\[2pt]
0 & l<k &
\end{array}\right.$$
and, similarly \bgl\label{2*} f x'^k y'^l =
\left\{\begin{array}{lll}
(k+f)(k-1+f)\ldots & (k+1-l+f)f x'^{k-l}&k \ge l\\[2pt]
0 & k<l &
\end{array}\right.\egl
\begin{lemma}\label{131}
\begin{enumerate}
\item[(a)] $f x'^k y'^l f= 0 \; {\rm if}\; k \ne l$
\item[(b)] $fx'^{n+1} y'^{n+1}
f =y'^n f
 = (n+1+f) fx'^n y'^n f$
\item[(c)] $fx'^n y'^n f = (n+f)(n-1+f)
\ldots (2+f)(1+f)f^2$
\end{enumerate}
\end{lemma}
\begin{proof}
Using the relations $[x',y'] = 1+f$ and $x'f = fy' = 0$ we
obtain
$$
x'^k y'^l = x'^{k-1} y' x' y'^{l-1} + x'^{k-1} y'^{l-1}
$$
$$
= x'^{k-1} y'^l x' + lx'^{k-1} y'^{l-1} + x'^{k-1} y'^{l-1} f
$$
$$
= x'^{k-1} y'^l x' + lx'^{k-1} y'^{l-1} + x'^{k-1} y'^{l-1} f
$$
Multiplying this by $f$ on both sides and using $x'f= 0$ gives
(a) and (b) (note also that $f$ commutes with $x'^n y'^n$ by
Lemma \ref{xy}. (c) follows from (b) by induction.\eproof By
Lemma \ref{w} every element of $I$ has a unique representation
of the form
$$
\sum y'^k z_{kl} x'^l$$ where $z_{kl} \in f \Cz[x'y']$. We write
$I_0$ for the ideal $f \Cz[x'y']$ in $\Cz[x'y']$ and consider
this algebra from now on as a subalgebra of $\mc C^\infty [0,1]$
by mapping $x'y'$ to the function $t$ and $f$ to the function
$t-1$ (this corresponds to the representation of $W'$ constructed
in section \ref{fund}).\\ We decompose $W'$ as $W' = I \oplus W$
using the linear splitting $s: W \ni y^nx^m \longmapsto y'^n(1+f)
x'^m\in W'$. Every element of $W'$ can then be uniquely
represented in the form $$\sum y'^k z_{kl} x'^l$$ with, this
time, $z_{kl} \in I_1$, where $I_1$ is the subalgebra
$(1+f)\Cz[x'y']$ of $\mc C^\infty [0,1]$ ($I_0$ consist in this
identification of the polynomial functions on $[0,1]$, that
vanish at 1). Here, the elements of the direct summand $W$
correspond to sums where all $z_{kl}$ are equal to $1+f =
x'y'-y'x'$.\\
Let $\beta$ be a submultiplicative seminorm on $I_1$ such
that $\beta(1) = 1$ and $\beta(1+f) \le 2$ (this is the case
for the restrictions of the canonical submultiplicative
seminorms on $\mc C^\infty [0,1]$ to $I_1$). We consider
seminorms $\beta_{\varphi}$ defined, using such a seminorm
$\beta$ and a monotone increasing function $\varphi: \Nz
\rightarrow \Rz_+$, by
$$
\beta_{\varphi}\big(\sum y'^k z_{kl} x'^l\big) = \sum
\varphi(k) \varphi(l) \beta(z_{kl})
$$
Let now $a,b \in I$. They can be written uniquely as $a =
\sum y'^k fz_{kl}x'^l$,\\ $b = \sum y'^i fw_{ij} x'^j$ with
$z_{kl}, w_{ij} \in I_1$. We have
$$
y'^k fz_{kl} x'^l y'^i fw_{ij} x'^j = \left\{\begin{array}{ll}
y'^k z_{kl}w_{lj} f^2(l+f)\ldots (2+f)(1+f) x'^j &l = i\\[2mm]
0  & l \ne i
\end{array}\right.$$
It follows that $$\begin{array}{l} \beta_{\varphi}(y'^k
fz_{kl} x'^l y'^i fw_{ij} x'^j) = \varphi(k) \varphi(j)
\beta(z_{kl} w_{lj} f^2(l+f) \ldots (1+f))\\[2mm]
\le \varphi(k) \beta(z_{kl}) \, 4(l+2) \ldots (1+2)
\beta(w_{lj}) \varphi(j) \le \varphi(k) \beta(z_{kl}) 2(l+2)!
\beta(w_{lj}) \varphi(j)\end{array}$$ if $l=i$ (and 0
otherwise). This last expression can be estimated by
$$
\varphi(k)\beta(z_{kl}) \varphi(l) \varphi(i) \beta(w_{ij})
\varphi(j)
$$
if $\varphi(l)^2 \ge 2(l+2)!$.\\
We thus obtain:   If $\varphi(l) \ge \sqrt{2(l+2)!}$, then
$\beta_{\varphi}$ is submultiplicative on $I$.\mn Let now $a=
\sum y'^k fz_{kl} x'^l \in I$ and $b= \sum y'^i (1+f) x'^j
\in s(W)$. We have from equation (\ref{2*}) at the beginning
of this section
$$
y'^k fz_{kl}x'^l y'^i(1+f) x'^j = \left\{\begin{array}{ll}
y'^k z_{kl} f(l+f) \ldots (l+1 - i + f)
x'^{l-i} (1+f)x'^j&l \ge i\\[2mm]
0 & l<i
\end{array}\right.
$$
whence
$$
\beta_{\varphi}(y'^k z_{kl} fx'^l y'^i (1+f) x'^j) \le
\varphi(k) \beta(z_{kl}) 4(l+1) \ldots (l+2-i) \varphi(l-i+j)
$$
(provided that $l \ge i$).\\
This, in turn, can be estimated by
$$
\varphi(k)\beta(z_{kl}) 4(l+1)! \varphi(2l)\psi(j)
$$
if $\psi(j) \ge \sup_{l \in \Nz}
\frac{\varphi(l+j)}{\varphi(2l)}$.\mn Setting $\varphi'(l) =
4(l+1)! \varphi(2l)$ this takes the form $\varphi(k)
\beta(z_{kl}) \varphi'(l) \psi(j).$ There is an analogous
estimate for products of the form $y'^i(1+f) x'^j y'^k z_{kl}
fx'^l$. As a consequence we got putting
$$
\alpha_{\psi}\Big(\sum_{i,j}y'^i(1+f) x'^j\Big) = \sum \psi(i)
\psi(j)
$$
with $\psi$ as above, the inequalities
$$
\beta_{\varphi}(ab) \le \beta_{\varphi'}(a) \alpha_{\psi}(b)
$$
$$
\beta_{\varphi}(ba) \le \alpha_{\psi}(b) \beta_{\varphi'}(a)
$$
for $a \in I,\, b \in s(W)$.\mn Let now finally $b_1 = s(c_1),
b_2 = s(c_2) \in s(W)$.\\ Then $b_1b_2 = s(c_1c_2) +
\omega(c_1, c_2)$, where $\omega (a,b)=s(ab)-s(a)s(b)$. We can
choose $\psi'$ such that
$$
\alpha_{\psi}(s(c_1c_2)) + \beta_{\varphi}(\omega(c_1,c_2))
\le \alpha_{\psi'}(b_1) \alpha_{\psi'}(b_2)
$$
Let then $\alpha, \alpha'$ on $W'= I + s(W)$ be given by\\
$$\begin{array}{lll}
\alpha(a+b) & = & \beta_{\varphi}(a) + \alpha_{\psi}(b)\\
\alpha'(a+b) & = & \beta_{\varphi'}(a) + \alpha_{\psi'}(b)
\end{array}$$\\
It follows that $\alpha(zw) \le \alpha'(z)\alpha'(w)$ for all
$z,w \in W'$. Note also that, for suitable $\psi$,
$\alpha_\psi$ dominates any given seminorm on $s(W)$.\mn In
conclusion we have proved the following theorem
\btheo\label{topW'} There is a locally convex structure on
$W'$ given by a family of seminorms $(\alpha_{\varphi,\psi})$
satisfying the following conditions:
\begin{itemize} \item the multiplication on $W'$ is
continuous\item as a locally convex space $W'$ is the direct sum
of $I$ and $W$ where $W$ carries the fine locally convex topology
and the restrictions of all seminorms to $I$ are
submultiplicative\item the locally convex structure on $I_1$
(defined as after \ref{131}) inherited from $W'$ is the usual
Fr\'{e}chet topology on $\mc C^\infty [0,1]\supset
I_1$.\end{itemize} The completion $\bar{W}'$ of $W'$ gives rise
to a linearly split extension of locally convex algebras $$0\to
\bar{I}\to \bar{W}'\to W\to 0$$ where $\bar{I}$ is an
$m$-algebra. Moreover, the completion of $I_1\subset W'$ in this
topology is isomorphic to $\mc C^\infty [0,1]$. \etheo \bremark
To define the topology with these properties on $W'$, it suffices
to take on $I$ only seminorms $\beta_\varphi$ where $\varphi$ is
obtained from $\varphi_0(k)\defeq \sqrt {2(k+1)!}$ by iteration
of the operation $\varphi '(k)=4(k+1)!\varphi (2k)$.\eremark
\section{Crossed products}
Many computations of the $K$-theoretic invariants for
$C^*$-algebras carry over to the locally convex setting. In
this section we discuss the case of crossed products of
locally convex algebras by $\Zz$. Let $A$ be a unital locally
convex algebra (the non-unital case can be treated by
adjoining a unit) and $\alpha$ a continuous automorphism of
$A$. The crossed product $A\rtimes_\alpha \Zz$ is defined as
the complete unital locally convex algebra generated by $A$
together with an invertible element $u$ satisfying
$uxu^{-1}=\alpha (x)$ for all $x$ in $A$ (thus for the
topology of $A\rtimes_\alpha \Zz$ we allow all seminorms on
the algebra generated algebraically by $A$ and $u,u^{-1}$ for
which the map $A\ni a\mapsto au^n$ is continuous for each
$n\in \Zz$). \\To describe the $kk^{\rm alg}$-invariants for
$A\rtimes_\alpha \Zz$ we follow the argument outlined in
\cite{CuK}. We denote by $\mc T_\alpha$ the closed subalgebra
of $\mc T\hot A\rtimes_\alpha \Zz$ generated by $1\otimes A$
and $v\otimes u$, $v^*\otimes u^{-1}$. Since the subalgebra of
$\mc C^\infty (S^1)\hot ( A\rtimes_\alpha \Zz)$ generated by
$1\otimes A$ and $z\otimes u$, $z^{-1}\otimes u^{-1}$ is
isomorphic to $A\rtimes_\alpha \Zz$ ($z$ the canonical
generator of $\mc C^\infty (S^1)$), we obtain, by
restriction, a continuous homomorphism $\mc T_\alpha\to
A\rtimes_\alpha \Zz$ whose kernel is obviously dense in the
ideal $\mc K\hot A$ of $\mc T_\alpha$. Thus, we have a
linearly split extension
$$0\to \mc K\hot A \to \mc T_\alpha\to A\rtimes_\alpha \Zz\to
0$$ There is a natural inclusion map $j:A\to \mc T_\alpha$
defined by $j(x)=1\otimes x$ and a quasihomomorphism $(\id,
\mbox{Ad}(v\otimes 1))$ from $\mc T_\alpha$ to $\mc K\hot A$
defining elements $kk(j)$ in $kk^{\rm alg}_0(A,\mc T_\alpha)$
and $\tau = kk(\id,\mbox{Ad}(v\otimes 1))$ in $kk^{\rm
alg}_0(\mc T_\alpha , \mc K \hot A)\cong kk^{\rm alg}_0(\mc
T_\alpha , A)$. Here, $\mbox{Ad}(v\otimes 1)$ is defined by
$\mbox{Ad}(v\otimes 1)(x)= (v\otimes 1)x(v^*\otimes
1)$.\bprop The elements $kk(j)$ and $\tau$ are inverse to
each other.\eprop\bproof It is clear that $kk(j)\cdot \tau
=1_A$. To prove that $\tau\cdot kk(j)=1_{\mc T_\alpha}$ let
$\psi$ be as in the proof of \ref{pair}. Let $k$ be the
homomorphism from $\mc T_\alpha$ to $\mc T \hot\mc T_\alpha$
which is given by the restriction of $\psi\otimes\id
_{A\rtimes \Zz}$ to $\mc T _\alpha$. The product $\tau\cdot
kk(j)$ is given by $kk(k,\mbox{Ad}(v\otimes 1\otimes 1)k)$.
By \ref{pair}, the homomorphism $k: \mc T_\alpha \to \mc
T\hot \mc T_\alpha$ is diffotopic to the sum of Ad$(v\otimes
1\otimes 1)k$ and a homomorphism $k'$ defined by $k'(1\otimes
a)= e\otimes 1\otimes a$, for $a$ in $A$, and $k'(v\otimes
u)= e\otimes v\otimes u$. Therefore $\tau\cdot kk(j)= kk(k')=
1_{\mc T_\alpha}$. \eproof \bprop Let $\kappa$ be the
inclusion of $A$ into $\mc K\hot A$ and $s$ the inclusion of
$\mc K\hot A$ into $\mc T_\alpha$. Then the element
$kk(s\kappa)\cdot \tau$ of $kk^{\rm alg}(A,A)$ is equal to
$1_A - kk(\alpha^{-1})$.\eprop \bproof $j$ is the sum of
Ad$(v\otimes 1)$ and $s\kappa$. On $\mc T_\alpha$
$$\mbox{Ad}(v\otimes 1) = \mbox{Ad}(1\otimes u^{-1})\mbox{Ad}
(v\otimes u) = \hat{\alpha}^{-1}\mbox{Ad}(v\otimes u)$$ where
$\hat{\alpha}$ is defined on $\mc T_\alpha$ by
$\hat{\alpha}(1\otimes a) = 1\otimes \alpha (a)$ and
$\hat{\alpha}(v\otimes u) =v\otimes u$. Since Ad$(v\otimes
u)$ is inner on $\mc T_\alpha$, one has
$kk(\mbox{Ad}(v\otimes u))= 1_{\mc T_\alpha}$ (by a standard
argument, inner endomorphisms are diffotopic to id after
stabilization by $2\times 2$-matrices). Therefore
$kk(\mbox{Ad}(v\otimes 1))\cdot \tau =kk(\alpha^{-1})$ on $A$
and $1_A = kk(j)\cdot\tau =kk(s\kappa)\cdot \tau +
kk(\mbox{Ad}(v\otimes 1))\cdot \tau $.\eproof \btheo For
every locally convex algebra $D$ there is an exact sequence
$$\begin{array}{ccccc} kk^{\rm alg}_0( D, A) & \stackrel{\,\cdot
(1_A - kk(\alpha^{-1}))}{\lori} & kk^{\rm alg}_0( D,
 A) &\lori& kk^{\rm alg}_0(
D, A\rtimes_\alpha \Zz)
\\[3pt]
\uparrow  & & & & \downarrow
\\[2pt]
kk^{\rm alg}_1( D, A\rtimes_\alpha \Zz) & \longleftarrow &
kk^{\rm alg}_1( D, A) & \stackrel{\,\cdot (1_A -
kk(\alpha^{-1}))}{\longleftarrow} & kk^{\rm alg}_1( D,A)
\end{array}$$ and a similar exact sequence in the first
variable. \etheo \bproof This is the long exact sequence
associated with the extension $$0\to \mc K\hot A \to \mc
T_\alpha\to A\rtimes_\alpha \Zz\to 0$$ after identifying
$kk^{\rm alg}(D,\mc T_\alpha)$ with $kk^{\rm alg}(D,A)$ using
the previous two propositions.\eproof Consider the unital
algebra $A(D,u)$ generated by two elements $D$ and $u$, where
$u$ is invertible, satisfying the commutation relation
$[D,u]u^{-1}=1$. This relation is used by Connes as one of
the simplest examples of a spectral triple, see for instance
\cite{CoSp}. The relation can be rewritten as $uDu^{-1}=D-1$.
Thus $A(D,u)$ is the crossed product $\Cz[D] \rtimes_\alpha
\Zz$ where $\Zz$ acts on the polynomial ring $\Cz[D]$ by the
automorphism determined by $\alpha (D)=D-1$. Obviously
$kk(\alpha)=1$. Therefore the Pimsner-Voiculescu long exact
sequence gives for instance for the $K$-theory $kk^{\rm
alg}_0(\Cz,A(D,u))\cong R$ and $kk^{\rm alg}
_1(\Cz,A(D,u))\cong R$. By the way, this algebra is
isomorphic to the algebra generated by two elements $x$ and
$z$, where $z$ is invertible, satisfying the relation
$xz-zx=1$.
\\
There are many examples of algebras for which a computation of
the $kk^{\rm alg}$-invariants should be possible using
similar methods. Consider for instance the algebra $B(D,u)$
generated by two elements $u,D$, where $u$ and $D$ are both
invertible and satisfy the relation $[D,u]u^{-1}=1$. One
finds that $kk^{\rm alg}_0(\Cz,B(D,u))\cong R$ and $kk^{\rm
alg} _1
(\Cz,B(D,u))\cong R^2$.\\
Another crossed product construction that has been considered
for instance in \cite{Adler} is the crossed product by a
derivation. Assume that $A$ is a locally convex algebra and
$\delta$ a continuous derivation of $A$. Define
$A\ltimes_\delta\Nz$ as the algebra generated by $A$ and an
element $D$ satisfying the relation $[D,a]=\delta(a)$ for
each $a\in A$ equipped with the finest locally convex
topology such that the map $A\ni a\mapsto D^na\in
A\ltimes\delta$ is continuous for each $n$. In this way
$A\ltimes_\delta\Nz$ becomes a locally convex algebra. The
Weyl algebra is a special case of this construction with
$A=\Cz[x]$. It is an obvious conjecture that
$A\ltimes_\delta\Nz$ should be isomorphic to $A$ in $kk^{\rm
alg}$. The corresponding result in periodic cyclic homology is
known, see \cite{Khal}.
\section{$m$-algebra approximation for the $K$-theory
of locally convex algebras}\label{app} In this section we
define a variant of bivariant $K$-theory for locally convex
algebras which reduces to the theory developed in
\cite{CuDoc} when applied to $m$-algebras. In particular, in
this theory, we can control the coefficients and we have
$kk(\Cz,\Cz)=\Zz$. This
allows to obtain more precise information in many cases.\\
Notably we use it in combination with the results from the
previous sections to derive some new results for the Weyl
algebra which seem to be difficult to get differently.\\
Recall that we defined in \cite{CuDoc}, for a locally convex
vector space $V$, the tensor $m$-algebra $\hat{T}V$ (denoted
simply $TV$ in \cite{CuDoc}) as the completion of the
algebraic tensor algebra
$$
T V=V\,\oplus\, V\!\!\otimes\!\! V\,\oplus\, V^{\otimes 3}
\oplus\,\dots
$$
with respect to the family $\{\widehat{p}\}$ of seminorms,
which are given on this direct sum as $$
\widehat{p}=p\,+\,p\!\otimes\!p\,+\,p^{\otimes 3 }\,+\,\dots
$$
where $p$ runs through all continuous seminorms on $V$. The
seminorms $\widehat{p}$ are submultiplicative for the
multiplication on $TV$. The completion $\hat{T} V$ therefore
is an $m$-algebra.\\ If $B$ is an $m$-algebra, then the
natural homomorphism $TB\to B$ extends to a continuous
homomorphism $\hat{T}B\to B$ and we obtain the free extension
$$0\to \hat{J}B\to \hat{T}B\to B\to 0$$ in the category of $m$-
algebras. It comes equipped with a natural continuous linear
splitting $\sigma:B\to \hat{T}B$.\\ Let us now say that an
extension
$$0\to I \to E\mathop{\lori} \limits^{\mathop{\textstyle
\curvearrowleft}\limits^{\scriptstyle s}}  D  \to 0$$ of
arbitrary locally convex algebras equipped with a continuous
linear splitting $s$ is weakly admissible (or that $s$ is
weakly admissible), if, for any continuous homomorphism
$\varphi: B\to D$ from an $m$-algebra $B$, the natural
homomorphism $\tau : TB\to E$ given by
$b_1\otimes\ldots\otimes b_n \mapsto
s\varphi(b_1)s\varphi(b_2)\ldots s\varphi(b_n)$ extends to a
continuous homomorphism and thus to a morphism of extensions
\[
\begin{array}{ccccccccc} 0 &\to &  \hat{J}B & \to &
\hat{T}B  & \to &  B & \to & 0\\[0.1cm]
&      &  \downarrow & &
\quad\downarrow\tau & & \quad\downarrow \varphi & &\\[0.1cm]
0 & \to & I & \to & E & \to & D & \to & 0
\end{array}
\]
We say that $0\to I \to E\mathop{\lori}
\limits^{\mathop{\textstyle
\curvearrowleft}\limits^{\scriptstyle s}}  D  \to 0$ is
admissible (or that $s$ is admissible) if, for any
$m$-algebra $C$, the extension $$0\to I\hot C \to E\hot C
\mathop{\lori} \limits^{\mathop{\textstyle
\curvearrowleft}\limits^{\scriptstyle s\otimes\id_C}}  D\hot C
\to 0$$ is weakly admissible.
\begin{lemma}\label{schlange} Let $\sigma :B\to TB$ be the
canonical linear map and, given $x,y\in B$, let $\omega
(x,y)=\sigma (xy)-\sigma (x)\sigma (y)$. Let $\beta$ be a
submultiplicative seminorm on $B$ and $\alpha$ a seminorm on
$TB$ which is submultiplicative on $JB$ and satisfies
$$\alpha (\sigma (x))\leq\beta (x)\,,\quad \alpha (\omega (x,y))\leq
\beta (x)\beta (y)\,,\quad \alpha (\sigma (x)\omega
(y,z))\leq\beta (x)\beta (y)\beta (z)$$ for all $x,y,z$ in
$B$. Then $\alpha \leq \widehat{2\beta}$ for the canonical
submultiplicative seminorm $\widehat{2\beta}$ on $TB$
associated to the seminorm $2\beta$.\end{lemma} \bproof We
use the decomposition of $TB$ in terms of curvature forms for
the standard linear splitting $\sigma: B\to TB$.\\ Given an
algebra $B$, we denote by $\Omega B$ the universal algebra
generated by $x\in B$ with relations of $B$ and symbols
$dx,x\in B$, where $dx$ is linear in $x$ and satisfies
$d(xy)=xd(y) + d(x)y$. We do not impose $d1=0$, i.e., if $B$
has a unit, $d1\neq0$. $\Omega B$ is a direct sum of
subspaces $\Omega^n B$ generated by linear combinations of
$\; x_0dx_1 \,\dots\, dx_n\,$, and $\; dx_1 \,\dots\, dx_n\,
,\; x_j\in B$. This decomposition makes $\Omega B$ into a
graded algebra. As usual, we write $\deg(\omega)=n$ if
$\omega \in \Omega^n B$. \mn As a vector space, for $n\ge 1$,
\begin{equation}\label{iso1}
\Omega^nB \cong \widetilde B \otimes B^{\otimes n}\cong
B^{\otimes (n+1)} \oplus B^{\otimes n}
\end{equation}
(where $\widetilde B$ is $B$ with a unit adjoined, and
$1\otimes x_1\otimes\dots\otimes x_n$ corresponds to $
dx_1\dots dx_n$). The operator $d$ is defined on $\Omega B$ by
$$
\begin{array}{c}
d(x_0dx_1\dots dx_n)\, =\, dx_0dx_1\dots dx_n
\\[2mm]
d(dx_1\dots dx_n)\, =\, 0\hspace{1.85truecm}
\end{array}
$$
The connection between $TB$ and $\Omega B$ is based on the
following deformation
\begin{itemize}
\item[(1)] The following defines an associative product (
Fedosov product) on $\Omega B$:
\[
\omega_1 \circ \omega_2 \; = \; \omega_1 \omega_2 +
(-1)^{deg(\omega_1)} d \omega_1 d \omega_2
\] where $\omega_1 \omega_2$ is the ordinary product of $\omega_1, \omega_2$ in
$\Omega A$. Clearly, the even forms $\Omega^{ev}A$ form a
subalgebra of $(\Omega A, \circ)$ for $\circ$.
\item[(2)] Let $\sigma : A \to TA$ be the canonical linear inclusion of $A$ into the
tensor algebra over $A$ and, for $x,y \in A$, $\omega(x,y) =
\sigma (xy) - \sigma x\sigma y$. Then the map
\[
\alpha : x_0 dx_1 dx_2 \ldots dx_{2n-1} dx_{x_{2n}} \; \mapsto
\; \sigma (x_0) \omega (x_1,x_2) \ldots \omega (x_{2n-1},
x_{2n})
\] defines an \underline{isomorphism} of algebras $(\Omega^{ev}B, \circ)
\stackrel{\alpha}{\longrightarrow} TB$.
\item[(3)] $JB$ is exactly the image of the space of even forms
of degree $>0$ under the map in (3).
\end{itemize} Define now a seminorm $\bar{\beta}$
on $TB\cong \Omega^{ev}B$ by
$$\bar{\beta}=\beta_0 + \beta_2 + \beta_4 \ldots$$
where $\beta_{2n}$ is defined on $\Omega^{2n}B\cong
B^{\otimes 2n}\oplus B^{\otimes\, (2n+1)}$ by the projective
tensor powers of $\beta$. Thus, for instance
$$\bar{\beta} (\sigma (x_0) \omega (x_1,x_2) \ldots
\omega (x_{2n-1}, x_{2n})) = \beta (x_0)\beta (x_1)\ldots
\beta (x_{2n})$$ Since $\alpha$ is submultiplicative on
products of several $\omega$'s, it is clear that $\alpha \leq
\bar{\beta}$.\\ Let now $\omega=\omega (x_1,x_2)\ldots
\omega(x_{2n-1},x_{2n})$ be an element in $\Omega^{2n}$ that
does not start with a $\sigma (x)$. Then, for $b\in B$,
$\bar{\beta}(\sigma(b)\omega) =\beta (b)\bar{\beta}(\omega)$
and
$$\bar{\beta}(\sigma (b)\sigma (x_0)\omega =
\bar{\beta}(\omega(b,x_0)+\sigma (b x_0)\omega)\leq 2\beta
(b)\beta (x_0)\bar{\beta}(\omega)=2\beta (b)
\bar{\beta}(\sigma(x_0)\omega)$$ since $\beta$ is
submultiplicative. This shows that for arbitrary $\omega '$
in $TB$, we have $\bar{\beta}(\sigma (b)\omega ')\leq 2\beta
(b)\bar{\beta}(\omega ')$.\\ Thus, in particular, by induction
$$\bar{\beta}(\sigma (x_0)\sigma(x_1)\ldots
\sigma (x_n))\leq 2^n\beta (x_0)\beta (x_1)\ldots \beta
(x_n)$$\eproof \bremark Let $A$ be any locally convex algebra.
If we denote by $\bar{T}A$ the completion of the tensor
algebra $TA$ with respect to all seminorms $\bar{\beta}$,
defined as in the previous proof, for all continuous
seminorms $\beta$ on $A$, then $\bar{T}A$ is a locally convex
algebra and the extension
$$0\to \bar{J}A\to\bar{T}A\to A\to 0$$ is admissible.\eremark
\begin{lemma}\label{ideal} Assume that $0\to
I\to E\to D\to 0$ is an extension of locally convex algebras
with a continuous linear splitting. If $I$ is an $m$-algebra
then the extension is admissible (with the given or with any
other linear splitting).\\ The classifying map $\hat{J}D\to
I$ is unique up to homotopy in such an extension.
\end{lemma} \bproof Obviously it suffices to show that any such
extension is weakly admissible. If $\varphi :B\to D$ is a
continuous homomorphism from an $m$-algebra $B$ and $\mu$ is a
continuous seminorm on $E$ which is submultiplicative on $I$,
then the pulled back seminorm $\alpha =\mu\circ\tau$ in the
diagram
\[
\begin{array}{ccccccccc} 0 &\to &  JB & \to &
TB  & \to &  B & \to & 0\\[0.1cm]
&      &  \downarrow & &
\quad\downarrow\tau & & \quad\downarrow \varphi & &\\[0.1cm]
0 & \to & I & \to & E & \to & D & \to & 0
\end{array}
\] satisfies the conditions of Lemma \ref{schlange} with
respect to a continuous submultiplicative seminorm $\beta$ on
$B$. The assertion then follows from \ref{schlange}\eproof
Given a locally convex algebra $A$, we equip the tensor
algebra $TA$ with the coarsest locally convex topology that
makes all natural maps $\tau$ continuous in \bgl\label{un}
\begin{array}{ccccccccc} 0 &\to &  JA & \to &
TA  & \to &  A & \to & 0\\[0.1cm]
&      &  \downarrow & &
\quad\downarrow\tau & & \quad\downarrow \beta & &\\[0.1cm]
0 & \to & I & \to & E & \to & D & \to & 0
\end{array}
\egl where the second row is an arbitrary admissible extension
and $\beta :A\to D$ a continuous homomorphism between locally
convex algebras. We denote the completions of $TA$ and $JA$
with respect to this locally convex structure by $\hat{T}A$
and $\hat{J}A$. Thus the vertical arrows in (\ref{un}) extend
to continuous homomorphisms from the completions $\hat{T}A$
and $\hat{J}A$ and
$$0\to \hat{J}A\to \hat{T}A\to A\to 0$$ is a universal
admissible extension. Note however that, in general, the
classifying map $\hat{J}A\to I$ depends on the continuous
linear splitting.\\
The characteristic defining properties of $A\mapsto\hat{T}A$
are
\begin{itemize}
\item If $A$ is an $m$-algebra, then $\hat{T}A$ is the
$m$-algebra tensor algebra from \cite{CuDoc}. \item The map
$A\mapsto\hat{T}A$ is functorial. \item For every $m$-algebra
$C$, the natural map $T(A\otimes C)\to TA\otimes C$ extends to a
continuous map $\hat{T}(A\hot C)\to \hat{T}A\hot C$.
\end{itemize}
Obviously, if $B$ is an $m$-algebra, then, by universality,
$\hat{T}B$ coincides with the tensor $m$-algebra described
above and $\hat{T}B$ and $\hat{J}B$ are $m$-algebras. Thus, by
induction, $\hat{J}^nB$ is an $m$-algebra for all $n$.
\begin{definition} Let $ A$ and $B$ be locally convex
algebras. We define
$$ kk( A ,\,  B\,) =
\lim_{\mathop{\lori}\limits_{k}} \langle \hat{J}^k A ,\,
\mathcal K \hat{\otimes} B(0,1)^k\, \rangle
$$
\end{definition}
We also define, for $n$ in $\Nz$,
$$kk_n(A,B)=kk(\hat{J}^nA,B)\qquad kk_{-n}(A,B)=kk(A,B(0,1)^n)$$
The product $kk_n \times kk_m \to kk_{n+m}$ can be defined
exactly as for $kk^{\rm alg}$. Lemma \ref{well} carries over,
if the extension is admissible and $J$ is replaced by
$\hat{J}$, showing that the product is well-defined.\\ As
remarked above in \ref{gleich}, this definition of $kk$, even
though formally slightly different, reduces to the one given
in \cite{CuDoc} in the case where $A$ and $B$ are
$m$-algebras. In fact, there are obvious maps both ways
between the groups defined here and in \cite{CuDoc},
respectively. It is obvious that these maps are inverse to
each other. Thus, the theory $kk$ is an extension of the
theory defined in \cite{CuDoc} from the category of
$m$-algebras to the category of locally convex algebras. Note
that the definition of $\hat{J}A$ above is set up exactly in
such a way as to make the theory functorial for maps from an
$m$-algebra $B$ into an arbitrary locally convex algebra
$A$.\\ In particular, we have the following important and
nontrivial fact, cf. \cite{CuDoc}, 7.2. \btheo For every
Banach algebra $B$, $kk_*(\Cz,B)$ coincides with the usual
topological $K$-theory $K_*(B)$. In particular,
$kk_n(\Cz,\Cz)=\Zz$, for $n$ even, and $kk_n(\Cz,\Cz)=0$ for
$n$ odd.\etheo The arguments used above to prove the mapping
cone exact sequences and the resulting long exact sequences
carry over to prove the analogous statements for the theory
$kk$, restricting however to admissible extensions. The
crucial fact is that the extensions used in those arguments
are admissible. In particular, if $0\to I\to E\to B\to 0$ is
any admissible extension then $0\to I(0 ,1)\to E(0,1)\to
B(0,1)\to 0$ is admissible and the suspension extension $0\to
A(0 ,1)\to A[0,1)\to A\to 0$ is admissible for any locally
convex algebra $A$. We do not know if the Toeplitz extension
$0\to \mc K\hot A\to \mc T_0 \hot A\to A(0,1)\to 0$ is
admissible for arbitrary $A$. Therefore the proof of Bott
periodicity for $kk_n(A,B)$ carries over only if $A$ or $B$
is an $m$-algebra. Thus $kk_n(A,B)$ is Bott-periodic, if $A$
or $B$ is an $m$-algebra and we obtain
\begin{theorem}
Let $ D$ be an $m$-algebra. Every extension of locally convex
algebras
$$
E:\, 0\to  I\stackrel{i}{\lori}  A\, \mathop{\lori}\limits^q
 B\,\to 0
$$
admitting an \underline{admissible} linear section, induces
exact sequences in $kk( D,\cdot\,)$ and $kk(\,\cdot\,,  D)$
of the following form: \bgl \label{11exact}
\begin{array}{ccccc}
kk_0( D,  I) & \stackrel{\,\cdot kk(i)}{\lori} & kk_0( D,
 A) &\stackrel{\,\cdot kk(q)}{\lori} & kk_0(
D, B)
\\[3pt]
\uparrow  & & & & \downarrow
\\[2pt]
kk_1( D,  B) & \stackrel{\,\cdot kk(q)}{\longleftarrow} &
kk_1( D,  A) & \stackrel{\,\cdot kk(i)}{\longleftarrow} &
kk_1( D, I)
\end{array}
\egl and \bgl \label{22exact}
\begin{array}{ccccc}
kk_0( I, D) & \stackrel{kk(i)\cdot \,}{\longleftarrow} &
kk_0(  A,  D) & \stackrel{kk(q)\cdot \,}{\longleftarrow} &
kk_0( B, D)
\\[3pt]
\downarrow  & & & & \uparrow
\\[2pt]
kk_1( B,  D) & \stackrel{kk(q)\cdot\,}{\lori} & kk_1( A,
 D) &\stackrel{kk(i)\cdot\,}{\lori} & kk_1(
I,  D)
\end{array}
\egl The admissible linear section for $E$ defines a classifying map
$\hat{J} B \to  I$ and thus an element of $kk_1( B, I)$,
which we denote by $kk(E)$. The vertical arrows in
(\ref{11exact}) and (\ref{22exact}) are (up to a sign) given
by right and left multiplication, respectively, by this class
$kk(E)$.
\end{theorem}
The natural continuous homomorphism $JA\to \hat{J}A$ induces a
map $kk(A,B)\to kk^{\rm alg}(A,B)$. As a consequence of
\ref{Chern}, the map $\Zz =kk(\Cz,\Cz)\to kk^{\rm
alg}(\Cz,\Cz)=R$ is injective ($1$ is mapped to $1_A \in
HP_0(\Cz,\Cz)$ under the composition of this map and the Chern
character $ch$).\\ Let us now apply these results to the Weyl
algebra. Theorem \ref{topW'} shows that we can complete the
extension
$$0\to I\to W'\to W\to 0$$
studied in section \ref{fund} to a linearly split extension of
the form $$0\to \bar{I}\to \bar{W}'\to W\to 0$$ where $\bar{I}$
is an $m$-algebra. In fact, the same argument as for $I$ in
\ref{I} shows that $\bar{I}$ is Morita equivalent in the sense of
\ref{Morita} to the closure of $t(1-t)\Cz[t]$ in $\mc C^\infty
[0,1]$, which in turn is diffotopy equivalent to $\Cz (0,1)$.\\
The completed extension is admissible and we obtain a long exact
sequence in $kk$, which can be compared to the long exact
sequences in $kk^{\rm alg}$ and $HP_*$ in the following
commutative diagram
\[
\begin{array}{ccccccccccc} \ldots &\to &  kk_0(D,\bar{I}) &
\to & kk_0(D,\bar{W}')  & \to &
kk_0(D,W)&\stackrel{\delta}{\to} &kk_1(D,\bar{I})
& \to & \ldots\\[0.1cm]
&      &  \downarrow& &
\downarrow & & \downarrow & &\downarrow& &\\[0.1cm]
\ldots &\to &  kk^{\rm alg}_0(D,\bar{I}) & \to & kk^{\rm
alg}_0(D,\bar{W}') & \to & kk^{\rm alg}_0(D,W)&\to &kk^{\rm
alg}_1(D,\bar{I})
& \to & \ldots\\[0.1cm]&      &  \downarrow& &
\downarrow & & \downarrow & &\downarrow& &\\[0.1cm]
 \ldots &\to &  HP_0(D,\bar{I}) & \to &
HP_0(D,\bar{W}')  & \to &  HP_0(D,W)&\to &HP_1(D,\bar{I}) &
\to & \ldots
\end{array}
\]
For $D=\Cz$, we have $kk^{\rm
alg}_*(\Cz,\bar{W}')=HP_*(\Cz,\bar{W}')=0$ and
$kk_0(\Cz,\bar{I})=0$, $kk_1(\Cz,\bar{I})=\Zz$. We think that
it can be shown that $kk(\Cz,\bar{W}')=0$, but the proof
along the lines of the argument in section \ref{W'} would be
very technical, because the necessary continuity conditions
are hard to check. At any rate it is very easy to
see that the map $kk(\Cz,\bar{I})\to kk(\Cz,\bar{W}')$ is 0.\\
However the commutative diagram above contains already enough
information for some applications. It shows that
$kk_0(\Cz,W)=G\oplus \Zz$, where $G\cong
kk_0(\Cz,\bar{W}')=\Ker \delta$ (which is presumably 0) and
$\delta$ maps the $\Zz$-summand isomorphically onto $\Zz =
kk_1(\Cz,\bar{I})$. Moreover, the $\Zz$-summand is also the
image of $\Zz =kk_0(\Cz,\Cz)$ under the natural map $\Cz\to
W$ mapping a complex number to the corresponding multiple of
1. This follows from the commutative diagram
\[
\begin{array}{ccccccccc} 0 &\to &  \Cz(0,1) & \to &
\Cz(0,1]  & \to &  \Cz & \to & 0\\[0.1cm]
&      &  \downarrow & &
\downarrow & & \downarrow  & &\\[0.1cm]
0 & \to & \bar{I}& \to & \bar{W}' & \to & W & \to & 0
\end{array}
\]
(where the map $\Cz(0,1]\to W'$ maps the generator $t$ of
$\Cz(0,1]$ to $x'y'$ in $W'$) and the induced commutative
diagram in the long exact sequences for
$kk_*(\Cz,\,\cdot\,)$.\\
From this we can deduce for instance that the algebras
$M_k(W)$ and $M_l(W)$ of matrices over $W$ can not be
isomorphic if $k\neq l$. This could also be shown using ring
theoretic methods (I am indebted to J.Dixmier for pointing
this out to me) or from the fact that the algebraic $K$-group
$K_0(W)$ is $\Zz$ which follows from \cite{Quill}, \S 6,
Theorem 7, cf. also \cite{AlevLambre}. The result that we get
here is however stronger, since the invariants that we use are
diffotopy invariants. Thus, for instance, we can also show
that for any homomorphism $\varphi : M_k(\Cz)\to \mc K\hot W$
the image $\varphi (1_k)$ of the unit in $M_k(\Cz)$ can be
diffotopic to the unit $1_l$ of $M_l(W)\subset \mc K\hot W$,
only if $l$ is a multiple of $k$.\\
As another example it follows that $M_k(W\hot \mc O_\infty )$
can be isomorphic to  $M_l(W\hot \mc O_\infty )$ only if
$k=l$ (using the fact that $kk_1(\Cz,\bar{I}\hot\mc
O_\infty)=\Zz$), or that $W\hot \mc O_n$
can be isomorphic to $W\hot \mc O_n$ only if $n=m$,  for the
$C^*$-algebras $\mc O_\infty$ and $\mc O_n$ considered in
\cite{CuCMP}.

\end{document}